\documentclass[12pt,oneside]{amsart}

\usepackage[T1]{fontenc}
\usepackage[utf8]{inputenc}
\usepackage{lmodern}
\usepackage[margin=1in]{geometry}
\usepackage{amsmath,amssymb,amsthm,mathtools}
\usepackage{enumitem}
\usepackage{microtype}
\usepackage[colorlinks=true,linkcolor=blue,citecolor=blue,urlcolor=blue]{hyperref}
\usepackage{todonotes}

\usepackage{tabularx,environ}

\makeatletter
\newcommand{\problemtitle}[1]{\gdef\@problemtitle{#1}}
\newcommand{\probleminput}[1]{\gdef\@probleminput{#1}}
\newcommand{\problemquestion}[1]{\gdef\@problemquestion{#1}}
\NewEnviron{algproblem}{
  \problemtitle{}\probleminput{}\problemquestion{}
  \BODY
  \par\addvspace{.5\baselineskip}
  \noindent
  \begin{tabularx}{\textwidth}{@{\hspace{0em}} l X c}
    \hline\hline
    \multicolumn{2}{@{\hspace{0em}}l}{\@problemtitle} \\
    \textbf{In:} & \@probleminput \\%
    \textbf{Question:} & \@problemquestion\\
    \hline\hline
  \end{tabularx}
  \par\addvspace{.5\baselineskip}
}
\makeatother

\numberwithin{equation}{section}
\hypersetup{pdftitle={Quantum polynomial-time conjugacy in finitely generated metabelian groups},
 pdfsubject={Conjugacy, affine algebras, exponent lattices, quantum algorithms},
 pdfkeywords={metabelian groups, conjugacy, BQP, exponent lattices}}
\newtheorem{theorem}{Theorem}[section]
\newtheorem{proposition}[theorem]{Proposition}
\newtheorem{lemma}[theorem]{Lemma}
\newtheorem{corollary}[theorem]{Corollary}
\theoremstyle{definition}

\theoremstyle{remark}

\newcommand{\Z}{\mathbb Z}

\newcommand{\MN}{\mathbb N}
\newcommand{\Q}{\mathbb Q}
\newcommand{\F}{\mathbb F}

\DeclareMathOperator{\Ann}{Ann}
\DeclareMathOperator{\Rel}{Rel}

\DeclareMathOperator{\rank}{rank}
\DeclareMathOperator{\im}{im}

\DeclareMathOperator{\ord}{ord}
\DeclareMathOperator{\Spec}{Spec}
\DeclareMathOperator{\Min}{Min}
\DeclareMathOperator{\adj}{adj}

\DeclareMathOperator{\Frac}{Frac}

\DeclareMathOperator{\Alg}{Alg}
\DeclareMathOperator{\mono}{mono}

\newcommand{\gp}[1]{{\left\langle #1 \right\rangle}}
\newcommand{\gpr}[2]{{\left\langle #1 \mid #2 \right\rangle}}

\newcommand{\Set}[2]{\left\{\, #1 \;\middle|\; #2 \,\right\}}

\def\CM{{\mathcal M}}

\title[Quantum polynomial-time conjugacy in metabelian groups]{Quantum polynomial-time conjugacy\\in metabelian groups}
\author{Alexei Miasnikov and Alexander Ushakov}
\date{}
\address{Department of Mathematical Sciences, Stevens Institute of Technology, Hoboken NJ 07030}\email{amiasnik,aushakov@stevens.edu}

\begin{document}
\maketitle
\begin{center}
\small\itshape
Dedicated to the memory of\\
Vladimir N. Remeslennikov,\\
Vitaly A. Roman'kov, and\\
Gennady A. Noskov,\\
mathematicians of Dostoevsky Omsk State University
who made important contributions\\
to the algorithmic theory of metabelian groups.
\end{center}
\begin{abstract}
We prove that the conjugacy problem in finitely
presented metabelian groups with a fixed number
of generators is solvable in bounded-error quantum
polynomial time, uniformly in the defining relations
of the group. Both the group presentation and the
words representing the elements are part of the input.

The main ingredient is a quantum polynomial-time algorithm for
computing all multiplicative relations among a fixed number of units
in a commutative affine $\Z$-algebra presented with a fixed number
of variables. The algorithm reduces the problem to computations
in finite-dimensional algebras over $\Q$ and finite rings by
identifying the lattice of exponent vectors whose corresponding
products are algebraic over the coefficient field.
An effective arithmetic splitting combines the rational and
finite-ring computations.
The proof uses effective bounds for polynomial equations,
algorithms for finite-dimensional commutative algebras,
quantum integer factorization, and finite abelian
hidden-subgroup algorithms.

\end{abstract}

\section{Introduction}

The conjugacy problem is one of the classical algorithmic
problems in group theory. For finitely generated metabelian
groups, its decidability is well known, but the computational
complexity of the problem is less understood.
In this paper, we show that the conjugacy problem in
finitely presented metabelian groups with a fixed number
of generators can be solved in bounded-error quantum
polynomial time. The algorithm is uniform in the defining
relations of the group.

Our approach exploits the connection between metabelian
groups and modules over commutative rings.
If $G$ is metabelian, then its derived subgroup $A=G'$
is an abelian normal subgroup, and conjugation makes
$A$ a module over the commutative group ring
$\Z[G/A]$.
This allows us to reduce the conjugacy problem to an
orbit problem for elements of a finitely presented module,
and subsequently to a multiplicative subgroup membership
problem in a commutative affine $\Z$-algebra.

The main algebraic ingredient is a quantum polynomial-time
algorithm for computing multiplicative relations among
a fixed number of units in a finitely presented
commutative $\Z$-algebra with a fixed number of variables.
The algebra may have positive Krull dimension, nilpotent
elements, and integer torsion.
Our algorithm combines effective commutative algebra
with quantum algorithms for integer factorization
and finite abelian groups.

\subsection{Quantum polynomial-time algorithms in algebra}

Since Shor's quantum polynomial-time algorithms
for integer factorization and discrete logarithms
\cite{Shor}, quantum computation has been studied
extensively. Nevertheless, relatively few broad
classes of computational problems in algebra and
group theory are known to admit quantum
polynomial-time algorithms.
Many of the known results rely on reductions to
integer factorization, discrete logarithms,
or hidden subgroup problems.

In computational group theory, Watrous developed
quantum polynomial-time algorithms for order
computation, subgroup membership, and related
problems in finite solvable black-box groups
\cite{Watrous01}.
Cheung and Mosca gave a quantum polynomial-time
algorithm for decomposing finite abelian groups
\cite{CM}.
Babai, Beals, and Seress obtained randomized
polynomial-time algorithms for fundamental
problems in finite matrix groups over fields
of odd characteristic, using integer factorization
and discrete logarithm oracles \cite{BBS09}.
Together with Shor's algorithms, these yield
quantum polynomial-time algorithms for the
corresponding problems.

In algebraic number theory, Hallgren developed
quantum polynomial-time algorithms for Pell's
equation and the principal ideal problem in
real quadratic fields \cite{Hallgren02}.
Eisentr\"ager, Hallgren, Kitaev, and Song
extended quantum methods to computing unit
groups of number fields of arbitrary degree,
using a continuous hidden subgroup approach
\cite{EHKS14}.

Another notable example is the quantum
polynomial-time algorithm of Aharonov, Jones,
and Landau for additive approximation of
the Jones polynomial at roots of unity \cite{AJL09}.
For suitable roots of unity and approximation
parameters, this problem is
$\mathsf{BQP}$-complete.

\subsection{Main results and input conventions}
\label{se:results}

For a tuple $u=(u_1,\ldots,u_k)$ of units
in a commutative ring and $e\in\Z^k$, write
\[
u^e=u_1^{e_1}\cdots u_k^{e_k}.
\]
If $A$ is a commutative $K$-algebra, where
$K$ is a field, the \emph{algebraicity lattice}
of $u\in(A^\times)^k$ is
\begin{equation}\label{eq:Alg}
\Alg_K(u)=
\Set{e\in\Z^k}{u^e
\text{ is algebraic over }K}.
\end{equation}
For a commutative ring $A$, the
\emph{relation lattice} of $u\in(A^\times)^k$ is
\begin{equation}\label{eq:Rel}
\Rel_A(u)=
\Set{e\in\Z^k}{u^e=1}.
\end{equation}

Throughout the paper, polynomial running times
refer to bit complexity, with integer coefficients
and lattice coordinates encoded in binary.
The numbers of variables and units are fixed,
and the polynomial-time bounds may depend on
these numbers.
Polynomial degrees are measured numerically.
We allow dense polynomial representations,
as well as sparse representations whose numerical
degrees are polynomially bounded in the input length.
Our complexity bounds do not apply to unrestricted
sparse representations with binary-encoded exponents.

\begin{theorem}[Computing algebraicity lattices]
\label{thm:algebraicity-main}
Fix nonnegative integers $m,k$ and let
$K=\Q$ or $\F_p$.
Let
\[
A=K[X_1,\ldots,X_m]/(f_1,\ldots,f_r),
\]
and let $u=(u_1,\ldots,u_k)\in(A^\times)^k$
be specified by polynomial representatives
of the units and their inverses.
A basis of $\Alg_K(u)$ can be computed
in deterministic polynomial time.
Moreover, if $D\geq 2$ bounds the total degrees
of the input polynomials, then $\Alg_K(u)$
admits a basis $b_1,\ldots,b_s$ satisfying
\[
\|b_i\|_\infty
\leq k!\,D^{(m+1)k}
\qquad (1\leq i\leq s).
\]
\end{theorem}

\begin{theorem}[Computing relation lattices]
\label{thm:affine-main}
Fix nonnegative integers $m,k$.
Let
\[
A=\Z[X_1,\ldots,X_m]/(f_1,\ldots,f_r),
\]
and let $u=(u_1,\ldots,u_k)\in(A^\times)^k$
be specified by polynomial representatives
of the units and their inverses.
A basis of $\Rel_A(u)$ can be computed
in bounded-error quantum polynomial time.
\end{theorem}

\begin{corollary}[Subgroup membership and exponent representations]
\label{cor:affine-membership}
Fix nonnegative integers $m,k$.
Let $S$ and $u_1,\ldots,u_k$ be as in
Theorem~\ref{thm:affine-main}, and let $h\in S^\times$
be specified by polynomial representatives of
$h$ and $h^{-1}$.
There is a bounded-error quantum polynomial-time algorithm
that decides whether
\[
h\in\langle u_1,\ldots,u_k\rangle.
\]
If so, the algorithm computes a vector $e_0\in\Z^k$
such that
\[
\{e\in\Z^k:u^e=h\}
=e_0+\Rel_S(u_1,\ldots,u_k).
\]
\end{corollary}

\begin{theorem}\label{thm:group-main-uniform}
Fix $k\geq 1$. There is a bounded-error quantum
polynomial-time algorithm that, given a finite
presentation
\begin{equation}\label{eq:G}
G = \gpr{x_1,\ldots,x_k}{r_1,\ldots,r_m}_{\CM}
\end{equation}
of a metabelian group and two words $x,y$
over its generators, decides whether $x$
and $y$ are conjugate in $G$.
\end{theorem}

In Theorem~\ref{thm:group-main-uniform}, the number $k$
of generators is fixed, but the defining relations of $G$
are part of the input. The algorithm is uniform in these
relations and in the input words. Its polynomial running
time has degree depending only on $k$.

The reductions and the new lattice argument are proved below.
We use the standard polynomial algorithms for integer linear algebra
and Hermite and Smith normal forms, the effective algebraic bounds
specified in Section~\ref{sec:integer}, Lenstra--Silverberg's
Theorem~1.11(i) for explicitly represented finite-dimensional
commutative $\Q$-algebras \cite{LS}, and quantum factorization,
order finding and the finite abelian hidden-subgroup algorithm
\cite{Shor,CM}. These underlying algorithmic theorems are used with their stated
input representations; their proofs are not reproduced here.

We use right conjugation $g^h=h^{-1}gh$ and commutators
$[g,h]=g^{-1}h^{-1}gh$. Abelian normal subgroups are
written additively when regarded as modules.

\subsection{Related results}

Noskov established decidability of the conjugacy problem
in finitely generated metabelian groups \cite{Noskov}.
The reduction of conjugacy to orbit problems in modules
over integral group rings is part of this classical theory.
We give an explicit formulation of this reduction and
show that the required module presentations and
subsequent cyclic-module computations can be constructed
in polynomial time, uniformly over presentations with
a fixed number of generators.

Stronger classical complexity bounds are known for
several important classes of solvable groups.
Matthews studied conjugacy in free metabelian groups
using wreath products \cite{Matthews}.
Vassileva established polynomial-time algorithms for
conjugacy decision and search in free solvable groups
\cite{Vassileva}.
Ushakov subsequently obtained a quasi-quartic-time
Monte Carlo algorithm for conjugacy in free solvable
groups, with a complexity bound uniform in the
derived length, which need not be fixed
\cite{Ushakov14}.

These results were extended to a broader family of
groups by Gul, Sohrabi, and Ushakov
\cite{Gul-Sohrabi-Ushakov:2015}.
They studied groups of the form $F/N^{(d)}$,
where $F$ is a finitely generated free group and
$N$ is a normal subgroup of $F$.
Assuming that the word problem in $F/N$ is solvable
in polynomial time, they obtained polynomial-time
algorithms for the word, power, and conjugacy
problems in $F/N^{(d)}$, for every $d\geq 2$,
with a polynomial bound independent of $d$.
Free solvable groups arise as the special case $N=F$.

Miasnikov, Vassileva, and Wei\ss{} proved that
conjugacy decision in fixed free solvable groups
and fixed iterated wreath products of finitely
generated abelian groups belongs to DLOGTIME-uniform
$\mathsf{TC}^0$ \cite[Corollary 5.6]{MVW}.
Here the rank, derived length, and iterated wreath
construction are fixed, and the inputs are ordinary
words. Thus free metabelian groups admit stronger
classical complexity bounds, while our result applies
uniformly to arbitrary finitely presented metabelian
groups with a fixed number of generators.

Nilpotent groups provide another point of comparison.
Myasnikov and Wei\ss{} established uniform
$\mathsf{TC}^0$ algorithms for conjugacy and other
algorithmic problems in nilpotent groups of bounded
rank and class, using unary inputs and the
quotient-presentation model specified in their paper
\cite[Theorem 13]{MW}.
Their algorithms also produce conjugators.
For binary-exponent inputs, the corresponding
conjugacy computations are $\mathsf{TC}^0$-Turing
reducible to binary extended gcd.
In particular, conjugacy in a fixed finitely
generated nilpotent group, with elements represented
by ordinary words, belongs to $\mathsf{TC}^0$.
Nilpotent groups of arbitrary class form a related
comparison class rather than a subclass of
metabelian groups.

For certain split metabelian groups of finite
Pr\"ufer rank, Gryak, Kahrobaei, and
Mart\'inez-P\'erez investigated the conjugacy
search problem \cite{GKM}.
They obtained an exponential-time algorithm in
the general setting, polynomial-time algorithms
for a subfamily, and a reduction to discrete
logarithms for another family.
The exponential bound is an upper bound for their
algorithm, not a complexity lower bound.
These results illustrate the role of multiplicative
arithmetic in conjugacy problems, although they
concern more restricted classes of groups.

The algebraic component of our approach is closely
related to algorithms for multiplicative relations.
Ge developed polynomial-time algorithms for
multiplicative relations among algebraic numbers
\cite{Ge}.
Lenstra and Silverberg established algorithms for
explicitly represented finite-dimensional
commutative $\Q$-algebras, including nonreduced
algebras \cite[Theorem 1.11]{LS}.
Their input consists of rational structure constants
relative to a vector-space basis.
We apply their result only after constructing an
appropriate finite-dimensional subalgebra with an
exact representation of polynomial size.

Kreuzer and Walsh developed algorithms for
computing binomial parts of polynomial ideals and
related exponent and unit lattices over fields,
including affine algebras \cite{KWbinomial}.
A unit lattice records exponents for which a product
belongs to the coefficient field, whereas an
exponent lattice records products equal to $1$.
Both notions differ from the auxiliary lattice of
exponents producing algebraic elements that we
use below.

In subsequent work, Kreuzer and Walsh also
computed exponent lattices in arbitrary finitely
generated $\Z$-algebras \cite[Section 3]{KWunit},
although the later sections of their paper focus
on algebras finite as $\Z$-modules.
Their approach separates rational and torsion
contributions and treats higher prime powers by
lifting.
Thus the underlying algebraic reductions are
already known.
Our contribution is to establish polynomial
bit-complexity bounds for these computations in
the affine presentation model, with a fixed number
of variables and a fixed number of specified units.

The quantum ingredients are standard algorithms
for integer factorization, order finding, and
the finite abelian hidden-subgroup problem
\cite{Shor,CM}.
We apply them to finite rings with explicit
canonical encodings.
Our results establish quantum upper bounds;
we do not claim classical lower bounds or that
all quantum steps reduce to integer factorization.

\subsection{Outline of the argument}
\label{se:outline}

The proof has two main parts. We first develop an algorithm
for computing multiplicative relations among units in
commutative affine $\Z$-algebras. We then reduce the conjugacy
problem in a finitely presented metabelian group 
with a fixed number of generators
to a multiplicative subgroup-membership problem in such
an algebra.

For the algebraic part, let
\[
S=\Z[X_1,\ldots,X_m]/I
\]
be a commutative affine $\Z$-algebra, and let
$u_1,\ldots,u_k\in S^\times$.
We assume that $m$ and $k$ are fixed and that the
units and their inverses are given by polynomial
representatives. Our goal is to compute a basis of
the relation lattice
\[
\Rel_S(u_1,\ldots,u_k)
=
\Set{e\in\Z^k}{u_1^{e_1}\cdots u_k^{e_k}=1}.
\]

Section~\ref{sec:integer} establishes the necessary
algorithmic preliminaries. These include polynomial-time
computations in modules over polynomial rings and group
rings, effective degree bounds, algebraicity tests,
and an arithmetic splitting for affine $\Z$-algebras.

In Section~\ref{sec:algebraic-directions}, we consider
an affine algebra over $K=\Q$ or $\F_p$ and the lattice
\[
\Alg_K(u_1,\ldots,u_k)=
\Set{e\in\Z^k}{u^e\text{ is algebraic over }K}.
\]
We prove that this lattice has a basis with polynomially
bounded coordinates. Since the number of units is fixed,
the lattice can be computed by enumerating polynomially
many exponent vectors and testing algebraicity.
Its basis vectors give units generating a
finite-dimensional algebra over $K$.
In the rational case, we construct this algebra
explicitly and compute its multiplicative relations
using the Lenstra--Silverberg algorithm.

Section~\ref{sec:primepower} treats affine algebras
over $\Z/p^a\Z$. We first compute the algebraicity
lattice after reduction modulo $p$ and lift suitable
polynomial relations to characteristic $p^a$.
This allows us to construct an exact finite subring
containing the relevant units.
Their multiplicative relations are then computed
using quantum order finding and the finite abelian
hidden-subgroup algorithm.

In Section~\ref{sec:integral-exp}, we combine the
rational and prime-power cases. An effective arithmetic
splitting gives an injective homomorphism
\[
S\longrightarrow (S\otimes\Q)\times S/cS
\]
for a computable positive integer $c$.
Quantum integer factorization and the Chinese remainder
theorem reduce the relation-lattice problem over $S$
to the previously treated cases.
The resulting lattices are intersected using
polynomial-time integer linear algebra.
This proves Theorem~\ref{thm:affine-main}.
Section~\ref{sec:affine-membership} derives the
corresponding multiplicative subgroup-membership
algorithm.

The group-theoretic reduction begins in
Section~\ref{sec:group}.
For a metabelian group \(G\) given by an input presentation 
\eqref{eq:G} with a fixed number of generators, we put
\[
A=G',\qquad Q=G/A,\qquad R=\Z Q.
\]
Using the $R$-module structure of $A$, we reduce
the conjugacy problem for two words to an orbit
membership problem
\[
f=e\cdot q,\qquad q\in Q,
\]
in an explicitly presented $R$-module.
The construction and all required module computations
take deterministic polynomial time.

In Section~\ref{se:affine}, we reduce this orbit
membership problem to multiplicative subgroup
membership
\[
h\in\langle u_1,\ldots,u_k\rangle
\]
in a commutative affine $\Z$-algebra.
Since the number $k$ of generators is fixed,
the number of polynomial variables
and the number of units are fixed.
The algorithm of Section~\ref{sec:affine-membership}
therefore applies, yielding the quantum polynomial-time
conjugacy algorithm in Section~\ref{se:final-proof}.

\section{Preliminaries: polynomial-time computations 
in polynomial rings and modules}
\label{sec:integer}

This section collects the computational tools used in
the subsequent reductions and algorithms. We establish
polynomial-time bounds for linear algebra over polynomial
rings in a fixed number of variables, module membership,
and annihilator computations. We also develop the arithmetic
splitting and elimination procedures needed for the
exponent-lattice algorithm.

Fix $s\geq0$ and let
\[
P=\Z[Z_1,\ldots,Z_s].
\]
Throughout this section, the number of variables
$s$ is fixed. Polynomials are represented in
expanded form, with integer coefficients encoded
in binary. All complexity bounds refer to this
representation.

Following \cite{Asch}, the \emph{height} of a polynomial $f\in P$
is defined by
\[
h(f)=
\max
\Set{\log |a|}{a\neq0\text{ is a coefficient of }f},
\]
with $h(0)=0$. More generally, for a finite collection
$f_1,\ldots,f_k\in P$, we put
\[
h(f_1,\ldots,f_k)
=
\max_{1\leq i\leq k} h(f_i).
\]

\subsection{Polynomial-time computations in $P$-modules}
\label{se:Ptime-P-modules}

In this subsection, we establish polynomial-time algorithms
for solving linear systems and computing kernels over $P$.
These algorithms rely on degree bounds from \cite{Asch}
and standard integer linear algebra.
We then apply them to submodule membership and
annihilator computations in finitely presented $P$-modules.

\begin{proposition}\label{prop:linear}
The following problems can be solved in deterministic polynomial time:
\begin{itemize}
\item[(a)]
deciding solvability of a linear system over $P$;
\item[(b)]
finding a solution of a linear system over $P$ when one exists;
\item[(c)]
computing a finite generating set for the kernel of a matrix over $P$.
\end{itemize}
\end{proposition}

\begin{proof}
Let $F\in P^{m\times r}$ and $b\in P^m$.
Let $d\geq1$ be an upper bound on the degrees of the entries of
$F$ and $b$, and let $h$ be an upper bound on their heights.
The matrix dimensions are regarded as part of the input.

Suppose first that the system
\[
Fv=b
\]
is solvable over $P$.
By \cite[Theorem~6.5]{Asch}, there is a
constant $C$, depending only on $s$, such
that the system has a solution satisfying
\[
\deg v_i\leq C(2md)^C(h+1)
\qquad (1\leq i\leq r).
\]
Since $s$ is fixed, this bound is polynomial
in $m$, $d$, and $h+1$.
Set
\[
D=\left\lceil C(2md)^C(h+1)\right\rceil.
\]
Write the coordinates of $v$ as polynomials of degree at most
$D$ with unknown integer coefficients. Since $s$ is fixed,
the number of monomials of degree at most $D$ in
$Z_1,\ldots,Z_s$ is
\[
\binom{D+s}{s}=O(D^s),
\]
and hence the total number of unknown coefficients is polynomial
in the input size. 
Comparing coefficients in $Fv=b$ therefore
yields an integer linear system
\[
Ax=c
\]
of polynomial size, whose coefficients have polynomial bit length.
Using, for example, Smith or Hermite normal form, one can decide
in polynomial time whether this system has an integral solution.
When a solution exists, these algorithms also produce one whose
coordinates have polynomial bit length. Substituting these
coordinates into $v$ gives a solution of $Fv=b$ of polynomial
size. This proves (a) and (b).

It remains to prove (c). By \cite[Corollary 5.4]{Asch}, the $P$-module
\[
\ker_P F=\{v\in P^r:Fv=0\}
\]
admits a generating set consisting of polynomial vectors of degree
at most
\[
(2md)^{2^{O(s^2)}}.
\]
Since $s$ is fixed, we may choose an
explicit degree bound $D_0$ that is
polynomial in $m$ and $d$.

Write a polynomial vector $v\in P^r$ of degree at most $D_0$
with unknown integer coefficients. Since $s$ is fixed and $D_0$
is polynomially bounded, comparing coefficients in
\[
Fv=0
\]
yields a homogeneous system of linear equations over $\Z$ of
polynomial size and polynomial coefficient bit length.
Computing a $\Z$-basis of its integer kernel in polynomial time
gives vectors $v_1,\ldots,v_\ell\in\ker_P F$ that span over $\Z$
all elements of $\ker_P F$ of degree at most $D_0$.

By the choice of $D_0$, the $P$-module $\ker_P F$ has a generating
set whose elements all have degree at most $D_0$. Each of these
generators therefore lies in the $\Z$-span of
$v_1,\ldots,v_\ell$. Hence
\[
\ker_P F=Pv_1+\cdots+Pv_\ell.
\]
The computation and its output have polynomial bit complexity.
\end{proof}

\begin{corollary}\label{cor:module-computations}
For finitely presented $P$-modules, the following problems can be
solved in deterministic polynomial time:
\begin{itemize}
\item[(a)]
deciding membership in a finitely generated submodule and,
when membership holds, computing a representation in terms of
its generators;
\item[(b)]
computing a finite generating set for the annihilator of
a given element.
\end{itemize}
\end{corollary}

\begin{proof}
Let $M$ be a finitely presented $P$-module
given by a presentation
\[
M=P^q/\im U,
\qquad U\in P^{q\times\ell}.
\]
The columns of $U$ give the defining relations, and elements of
$M$ are represented by vectors in $P^q$.
Suppose that the columns
of a matrix $V$ represent generators of a submodule of $M$.
An element represented by $b\in P^q$ belongs to this submodule
if and only if there exist vectors $c$ and $z$ over $P$ such that
\[
Vc+Uz=b.
\]
This is a linear system over $P$. By
Proposition~\ref{prop:linear}, its solvability can be decided in
polynomial time and, when it is solvable, a solution $c$ gives
the required representation of $b$.

Now let $a\in P^q$ represent an element $\bar a\in M$.
A polynomial $f\in P$ belongs to $\Ann_P(\bar a)$ if and only if
\[
fa\in\im U,
\]
or equivalently, if there exists $z$ such that
\[
[\,a\mid -U\,]
\binom{f}{z}=0.
\]
By Proposition~\ref{prop:linear}, generators of the kernel of this
matrix can be computed in polynomial time. Projecting these
generators onto their first coordinates gives a generating set
for $\Ann_P(\bar a)$.
\end{proof}

\subsection{Polynomial-time computations in $\Z Q$-modules}
\label{se:Ptime-ZQ-modules}

Fix $k\geq 0$, and let $Q$ be an abelian group
given by a finite presentation
\begin{equation}\label{eq:Q-presentation}
Q=\langle x_1,\ldots,x_k \mid r_1,\ldots,r_m\rangle_{\mathrm{Ab}}
\end{equation}
in the variety of abelian groups, where
$r_1,\ldots,r_m$ are group words over
$x_1,\ldots,x_k$.
Put
\[
R=\Z Q,
\qquad
P=\Z[X_1,Y_1,\ldots,X_k,Y_k].
\]
For a group word $w$ over $x_1,\ldots,x_k$,
let $\sigma_{x_i}(w)$ denote the exponent sum
of $x_i$ in $w$, and define the monomial
\[
\mono(w) =
\bigg(\prod_{\substack{1\le i\le k\\\sigma_{x_i}(w)>0}}
X_i^{\sigma_{x_i}(w)}\bigg)
\bigg(\prod_{\substack{1\le i\le k\\\sigma_{x_i}(w)<0}}
Y_i^{-\sigma_{x_i}(w)}\bigg).
\]
Since $Q$ is abelian, we have
\begin{equation}\label{eq:R-presentation}
R\cong P/J,
\qquad
J=\gp{X_iY_i-1\ (1\le i\le k),\
\mono(r_j)-1\ (1\le j\le m)},
\end{equation}
where $X_i$ and $Y_i$ represent
$x_i$ and $x_i^{-1}$, respectively.

Since the relators are given as group words,
the exponent sums $\sigma_{x_i}(r_j)$ can be
computed in polynomial time and satisfy
$|\sigma_{x_i}(r_j)|\le |r_j|$.
Thus the monomials $\mono(r_j)$ have
polynomial degree and expanded size.

\begin{corollary}\label{cor:group-ring}
For fixed $k$, the following problems can be solved
in deterministic polynomial time, uniformly in
the defining relators $r_1,\ldots,r_m$ of $Q$:
\begin{itemize}
\item[(a)]
deciding solvability of a linear system over $R$ and finding a
solution when one exists;
\item[(b)]
computing a finite generating set for the kernel of a matrix over $R$;
\item[(c)]
deciding membership in a finitely generated submodule and,
when membership holds, computing a representation in terms of
its generators;
\item[(d)]
computing a finite generating set for the annihilator of a
given element.
\end{itemize}
\end{corollary}

\begin{proof}
A finitely presented $R$-module can be viewed
as a finitely presented $P$-module by adjoining
the relations $ge_i=0$ for each generator
$g$ of $J$ and each module generator $e_i$.
This conversion takes polynomial time and
produces a presentation of polynomial size.
Linear systems and submodule membership over $R$
then become the corresponding problems over $P$.
Kernels over $R$ are computed by lifting the
matrices to $P$, adjoining generators of $J$
in each target coordinate, and projecting
the resulting $P$-module kernels to $R$.
Annihilators are computed similarly.

Thus all four problems reduce to the algorithms
of Proposition~\ref{prop:linear} and
Corollary~\ref{cor:module-computations}.
Since $P$ has the fixed number $2k$ of variables,
their running times are polynomial in the combined
input size, including the lengths of the
defining relators of $Q$.
\end{proof}

\begin{corollary}\label{cor:inverse-mod-ideal}
For fixed $k\ge 0$, there is a deterministic polynomial-time
algorithm that, given a finitely generated ideal
\[
I=(f_1,\ldots,f_\ell)
\]
of $R$ and an element $h\in R$, decides whether
$h+I$ is a unit in $R/I$ and, when it is,
computes a representative of its inverse.
The running time is polynomial in the combined
input size, including the defining relators of $Q$.
\end{corollary}

\begin{proof}
The element $h+I$ is a unit in $R/I$ if and only if
there exist $h',c_1,\ldots,c_\ell\in R$ such that
\[
hh'+\sum_{j=1}^{\ell}f_jc_j=1.
\]
This is a linear system over $R$.
By Corollary~\ref{cor:group-ring}, its solvability
can be decided in deterministic polynomial time,
uniformly in the presentation of $Q$, and a solution
can be computed whenever one exists.
In that case, $h'+I$ is the inverse of $h+I$
in $R/I$.
\end{proof}

\subsection{Bounded-degree ideal computations}
\label{se:ideal-bounded-degree}

For a polynomial ring $R$ over a field or
over $\Z$, let $R_{\leq d}$ denote the set
of polynomials of total degree at most $d$.

For computations over $\Q$, coefficients are
represented as reduced fractions with binary
numerators and denominators. Over $\F_p$,
coefficients are represented by integers modulo
$p$, with $p$ encoded in binary.

We begin with an auxiliary lemma that reduces
bounded-degree ideal computations to linear
algebra over the coefficient field.
It will be used both to compute rational
contractions of ideals and to test algebraicity.

\begin{lemma}\label{lem:bounded-intersection}
Let $K=\Q$ or $\F_p$, and let
\[
R=K[Z_1,\ldots,Z_s],
\qquad
J=(f_1,\ldots,f_r)\leq R,
\]
where $s$ is fixed.
Given an integer $E\geq0$ encoded in unary,
there is a deterministic polynomial-time algorithm
that computes a $K$-basis of
\[
J\cap R_{\leq E}.
\]
Over $\Q$, the computed basis has polynomial
coefficient bit length.
\end{lemma}

\begin{proof}
Let $D\geq2$ bound the degrees of the
generators of $J$, and put
$d=\max\{D,E\}$.
By Hermann's ideal-membership degree bound
\cite[Theorem~3.4]{Asch}, every polynomial
$h\in J$ of degree at most $E$ admits a
representation
\[
h=\sum_{\nu=1}^r a_\nu f_\nu,
\qquad
\deg a_\nu\leq (2d)^{2^s}.
\]
Set $H=(2d)^{2^s}+D$.
Since $s$ is fixed, both the multiplier-degree
bound and $H$ are polynomial in $d$.

Write $h$ and the multipliers $a_\nu$
with undetermined coefficients, subject
to the respective degree bounds $E$
and $(2d)^{2^s}$.
Comparing coefficients up to degree $H$
gives a homogeneous linear system
\[
Bu=Ca,
\]
where $u$ is the coefficient vector of $h$
and $a$ collects the coefficients of the
multipliers.
Here $B$ embeds the coefficients of $h$
into the coefficient space of polynomials
of degree at most $H$, while $C$ represents
the map
$(a_1,\ldots,a_r)\mapsto\sum_\nu a_\nu f_\nu$.

By Gaussian elimination, compute a matrix
$L$ whose rows form a basis of the left
nullspace of $C$. Then
\[
Bu\in\operatorname{im}C
\quad\Longleftrightarrow\quad
LBu=0.
\]
Thus $\ker(LB)$ is precisely the coefficient
space of
\[
J\cap R_{\leq E}.
\]
Since $s$ is fixed, the degree bounds imply
that all the linear systems have polynomial
dimensions and polynomial coefficient bit length.
Standard bounds for Gaussian elimination over
$\Q$ ensure that the matrices produced by
elimination and the resulting basis vectors
also have polynomial coefficient bit length.
Over $\F_p$, field operations have polynomial
bit complexity in $\log p$.
Thus the required basis can be computed in
deterministic polynomial time.
\end{proof}

\subsection{Arithmetic splitting}
\label{se:arithmetic-splitting}

We now establish an effective separation of
rational and modular information in a
finitely generated $\Z$-algebra.
Theorem~\ref{th:arithmetic} embeds $P/I$ into
the product of its rationalization and a modular
quotient, where the modulus has polynomial bit length.
Consequently, multiplicative relations in $P/I$,
including exponent-lattice relations, can be
recovered by computing them in these two quotients
and intersecting the resulting lattices.
A similar separation into rational and
finite-characteristic parts appears in the
exponent-lattice algorithm of Kreuzer and Walsh
\cite[Section~3]{KWunit}.

Let $I=(f_1,\ldots,f_r)\leq P$, and put
\[
I_\Q=I\Q[Z_1,\ldots,Z_s],
\qquad
J=I_\Q\cap P.
\]
The rationalization of $P/I$ is naturally
isomorphic to
\[
(P/I)\otimes_\Z\Q
\cong \Q[Z_1,\ldots,Z_s]/I_\Q.
\]
The kernel of the natural map
\[
P/I\longrightarrow
\Q[Z_1,\ldots,Z_s]/I_\Q
\]
is $J/I$. Thus $J/I$ is precisely the additive
$\Z$-torsion subgroup of $P/I$, and $P/J$ is
its torsion-free quotient.
We first compute generators of $J$ and then
find an integer that annihilates $J/I$.

\begin{lemma}[Rational contraction]
\label{lem:rational-contraction}
There is a deterministic polynomial-time
algorithm that computes generators of $J$.
\end{lemma}

\begin{proof}
Let $D\geq2$ bound the degrees of
$f_1,\ldots,f_r$.
By \cite[Proposition~4.7]{Asch}, there is an
explicit degree bound
\[
\tau\leq (2D)^{2^{O(s^2)}}
\]
such that $J$ is generated by polynomials
of degree at most $\tau$.
Since $s$ is fixed, $\tau$ is polynomially
bounded in $D$.

Consider the rational vector space
\[
V=I_\Q\cap
\Q[Z_1,\ldots,Z_s]_{\leq\tau}.
\]
By the proof of Lemma~\ref{lem:bounded-intersection},
we can compute in polynomial time a rational
matrix $M$ whose kernel, under the monomial
coefficient representation, is precisely $V$.

Clearing the denominators in each row of
$M$ produces an integer matrix $H$ with
\[
\ker_\Q H=\ker_\Q M=V.
\]
Under the corresponding identification
$P_{\leq\tau}\cong\Z^N$, we have
\[
V\cap P_{\leq\tau}
\cong\ker_\Z H.
\]
The matrix $H$ has polynomial dimensions
and polynomial coefficient bit length.
Using a polynomial-time integer-kernel
algorithm based on Hermite normal form
\cite{KB}, compute a $\Z$-basis
\[
j_1,\ldots,j_t
\]
of this module.
These polynomials generate $J$ as an ideal.
Indeed, by the choice of $\tau$, the ideal
$J$ has generators of degree at most $\tau$,
all of which belong to
$V\cap P_{\leq\tau}$.
\end{proof}

\begin{lemma}[Torsion annihilator]
\label{lem:torsion-annihilator}
There is a deterministic polynomial-time
algorithm that computes 
$c\in\MN$ of polynomial bit length such that
\[
cJ\subseteq I.
\]
\end{lemma}

\begin{proof}
By Lemma~\ref{lem:rational-contraction},
compute generators $j_1,\ldots,j_t$ of $J$.
If $J=0$, the conclusion holds with $c=1$.
Otherwise, for each $j_i$, find a representation
\[
j_i=\sum_{\nu=1}^r a_{i\nu}f_\nu,
\qquad
a_{i\nu}\in\Q[Z_1,\ldots,Z_s].
\]
Hermann's degree bound
\cite[Theorem~3.4]{Asch} gives a polynomial
bound on the degrees of the multipliers.
Consequently, these representations can
be found by solving polynomial-size
rational linear systems.
The resulting coefficients have polynomial
bit length, and their total number is
polynomially bounded.
Let $c>0$ be the least common multiple
of all their denominators.
Then $c$ has polynomial bit length and
\[
cj_i\in I
\qquad (1\leq i\leq t).
\]
Since the $j_i$ generate $J$, it follows that
$cJ\subseteq I$.
\end{proof}

\begin{theorem}[Arithmetic splitting]
\label{th:arithmetic}
There is a deterministic polynomial-time
algorithm that computes 
$c\in\MN$ of polynomial bit length such that
the natural homomorphism
\[
P/I\longrightarrow
(P/I)\otimes_\Z\Q
\times
(P/I)/c(P/I)
\]
is injective.
\end{theorem}

\begin{proof}
Put $A=P/I$ and let $T=J/I$ be its
additive $\Z$-torsion subgroup.
By Lemma~\ref{lem:torsion-annihilator},
compute a positive integer $c$ of polynomial
bit length satisfying $cJ\subseteq I$.
It follows that
\[
cT=0.
\]
The kernel of the natural homomorphism
\[
A\longrightarrow
A\otimes_\Z\Q\times A/cA
\]
is $T\cap cA$.
Let $z\in T\cap cA$ and write $z=cy$.
Since $z$ is torsion, there is a nonzero
integer $n$ such that
\[
0=nz=ncy.
\]
Hence $y\in T$, so $z=cy=0$.
Therefore the homomorphism is injective.
\end{proof}

\subsection{Algebraicity test}
\label{se:algebraicity-test}

We next give an effective test for algebraicity
in finitely generated algebras over $\Q$ and
$\F_p$. For a fixed number of variables,
the test runs in deterministic polynomial time
and, when an element is algebraic, produces
a monic annihilating polynomial with
polynomially bounded degree and coefficient
bit length.

\begin{lemma}
\label{lem:algebraicity-test}
Let $K=\Q$ or $\F_p$, and let
\[
A=K[Z_1,\ldots,Z_s]/I,
\]
where $s$ is fixed.
Given a polynomial representative of $v\in A$,
there is a deterministic polynomial-time
algorithm that decides whether $v$ is
algebraic over $K$.
If it is, the algorithm computes a monic
annihilating polynomial of polynomial degree
and coefficient bit length.
\end{lemma}

\begin{proof}
Let $v$ be represented by
$g\in K[Z_1,\ldots,Z_s]$.
Introduce a new variable $Y$ and consider
the graph ideal
\[
J=I+(Y-g)
\leq K[Z_1,\ldots,Z_s,Y].
\]
The evaluation homomorphism
\[
\operatorname{ev}_v:
K[Z_1,\ldots,Z_s,Y]\longrightarrow A,
\]
defined by
\[
Z_i\longmapsto Z_i+I,
\qquad
Y\longmapsto v,
\]
has kernel $J$. Its restriction to $K[Y]$
therefore has kernel $J\cap K[Y]$.
Consequently,
\[
v\text{ is algebraic over }K
\quad\Longleftrightarrow\quad
J\cap K[Y]\ne0.
\]

Let $D\geq2$ bound the degrees of the
generators of $J$.
By Dub\'e's degree bound
\cite[Corollary~8.3]{Dube}, applied with
an elimination ordering, if
$J\cap K[Y]\ne0$, then it contains
a nonzero polynomial of degree at most
\[
E=
\left\lceil
2\left(\frac{D^2}{2}+D\right)^{2^s}
\right\rceil.
\]
Since $s$ is fixed, $E$ is polynomially
bounded in $D$.

By Lemma~\ref{lem:bounded-intersection},
applied to the polynomial ring
$K[Z_1,\ldots,Z_s,Y]$ in $s+1$ variables,
we can compute a basis of
\[
V=
J\cap K[Z_1,\ldots,Z_s,Y]_{\leq E}
\]
in deterministic polynomial time.
Represent the elements of $V$ by their
coefficient vectors in the monomial basis.
Imposing the linear conditions that all
coefficients of monomials involving any
$Z_i$ vanish, compute a basis of
\[
V\cap K[Y]
=
J\cap K[Y]_{\leq E}.
\]
By the choice of $E$, the element $v$
is algebraic over $K$ if and only if
this subspace is nonzero.
In that case, choose a nonzero polynomial
$h(Y)$ in this subspace and divide by
its leading coefficient.
The resulting monic polynomial annihilates $v$.

Since $s$ is fixed, all degree bounds and
linear-system dimensions are polynomially
bounded. The coefficient bit-length bounds
follow from Lemma~\ref{lem:bounded-intersection}
and standard bounds for rational linear algebra.
Thus the algorithm runs in deterministic
polynomial time.
\end{proof}

\section{Relation lattices over fields}
\label{sec:algebraic-directions}

The goal of this section is to compute multiplicative relation lattices 
over \(\mathbb Q\) and \(\mathbb F_p\). 
Our main tool is the algebraicity lattice of the given units.
The key point is that the resulting lattice admits a basis whose
entries are bounded polynomially in the numerical degrees of the
input polynomials. Since the number of units is fixed, this permits
an exhaustive search over a polynomial-size set of exponent vectors.

Throughout this section, the number $m$ of polynomial variables
and the number $k$ of specified units are fixed. We consider an
affine $K$-algebra
\[
A=K[X_1,\ldots,X_m]/I,
\qquad
I=(f_1,\ldots,f_r),
\]
where $K=\Q$ or $K=\F_p$.
The polynomials $f_1,\ldots,f_r$, together with polynomial
representatives of the units $u_1,\ldots,u_k\in A^\times$
and their inverses, are given as input.
Let $D\geq2$ bound the total degrees of all these polynomials.
Over $\Q$, coefficients are binary rationals; over $\F_p$,
their encoding uses $O(\log p)$ bits.

\subsection{Algebraicity lattices in affine algebras}
    
Recall the algebraicity lattice $\Alg_K(u)$
defined in \eqref{eq:Alg}.
It is a saturated sublattice of $\Z^k$.
Indeed, products and inverses of algebraic
units are algebraic, and if $(u^e)^n$
is algebraic over $K$ for some $n\ne0$,
then so is $u^e$.
Put
\[
V=D^{m+1},
\qquad
B_0=k!V^k.
\]

\begin{proposition}\label{pr:alg-lattice}
The lattice $\Alg_K(u)$ has a basis consisting of vectors $b$
satisfying
\[
\|b\|_\infty\leq B_0.
\]
\end{proposition}

\begin{proof}
Let $\mathfrak p$ be a minimal prime ideal of $A$. Then
$A/\mathfrak p$ is an integral domain. We write
\[
X=\Spec(A/\mathfrak p)
\]
for the corresponding irreducible component of $\Spec(A)$.
Its coordinate ring is $A/\mathfrak p$, and its function field is
\[
K(X)=\Frac(A/\mathfrak p).
\]
Thus the images of the units $u_1,\ldots,u_k$ in
$A/\mathfrak p$ may be viewed as nonzero rational functions on
$X$.

We pass to a normal projective model in order to record all zeros
and poles of these rational functions, including those at infinity.
Let
\[
\overline X\subseteq\mathbb P^m_K
\]
be the projective closure of $X$, and let
\[
\nu:Y\longrightarrow\overline X
\]
be its normalization. Then $Y$ is an integral normal projective
variety and $\nu$ is finite and birational; see
\cite[Tag 0BXR]{Stacks}. In particular,
\[
K(Y)\simeq K(\overline X)\simeq K(X),
\]
so we may regard the $u_i$ as rational functions on $Y$.

Since $X$ is an irreducible component of $\Spec(A)$
and $I$ is generated by polynomials of degree at most $D$,
the affine B{\'e}zout inequality
\cite[Theorem~1]{Heintz83} gives
\[
\deg X\leq D^m.
\]
Here the degree of an affine variety is the degree
of its projective closure under the standard embedding
$\mathbb A^m_K\hookrightarrow\mathbb P^m_K$.
Consequently,
\[
\deg\overline X=\deg X\leq D^m.
\]

Suppose $s=\dim X\geq1$. Let $F_i\in K[X_1,\ldots,X_m]$
be a polynomial representative of $u_i$ of degree at most $D$,
and let $F_i^h$ denote its homogenization to degree $D$.
On the affine chart $X_0\neq0$, the polynomial $F_i$ is represented
by the homogeneous quotient
\[
\frac{F_i^h}{X_0^D}.
\]
Thus, after pulling back along
$\nu:Y\to\overline X$, the corresponding rational function on
$Y$ is
\[
u_i
=
\nu^*\left(\frac{F_i^h}{X_0^D}\right)
=
\frac{\nu^*F_i^h}{\nu^*X_0^D}.
\]

Put
\[
H=\nu^*\mathcal O_{\overline X}(1).
\]
Since $\nu$ is finite, $H$ is ample.
The sections $\nu^*F_i^h$ and $\nu^*X_0^D$ are nonzero sections
of $H^{\otimes D}$, and hence their zero sets define effective
Cartier divisors on $Y$. We control their multiplicities using
their degrees with respect to $H$:
\[
c_1(H)^{s-1}\cdot\operatorname{div}(\nu^*F_i^h)
=
c_1(H)^{s-1}\cdot\operatorname{div}(\nu^*X_0^D)
=
D\deg\overline{X}.
\]
Since $\nu$ is finite and birational, $\deg(\nu)=1$, and
\cite[Tag~0BEX]{Stacks} gives
$\deg_{H}(Y)=\deg_{\mathcal{O}(1)}(\overline{X})=\deg\overline{X}$.
Since $H$ is ample, every prime divisor $E$ on $Y$ has $H$-degree at least $1$.
The multiplicities of the prime divisors in $\operatorname{div}(\nu^* F_i^h)$ are non-negative integers whose $H$-degree-weighted sum equals $D \deg \overline{X}$. Since each weight is at least $1$, each individual multiplicity is at most $D \deg \overline{X}$, and the same holds for $\operatorname{div}(\nu^* X_0^D)$.
Since
\[
\ord_E(u_i)
=
\ord_E(\nu^*F_i^h)
-
\ord_E(\nu^*X_0^D),
\]
with both terms on the right between $0$ and
$D\deg\overline X$, we obtain
\begin{equation}\label{eq:valuation-bound}
\bigl|\ord_E(u_i)\bigr|
\leq D\deg\overline X
\leq D^{m+1}=V.
\end{equation}

For
\[
u^e=u_1^{e_1}\cdots u_k^{e_k},
\]
the additivity of $\ord_E$ gives
\[
\ord_E(u^e)
=
\sum_{i=1}^k e_i\ord_E(u_i).
\]
Thus the condition that $u^e$ have neither a zero nor a pole
along $E$ is a homogeneous integer linear equation in $e$.

For each positive-dimensional component $X$ corresponding to a
minimal prime $\mathfrak p$ of $A$, and each prime divisor $E$ on
the normalization of its projective closure, consider the row
\[
\bigl(\ord_E(u_1),\ldots,\ord_E(u_k)\bigr)\in\Z^k.
\]
Only finitely many of these rows are nonzero, since each $u_i$
has only finitely many zeros and poles. Let $M$ be the integer
matrix formed by all nonzero rows. By
\eqref{eq:valuation-bound}, every entry of $M$ has absolute value
at most $V$. We claim that
\[
\Alg_K(u)=\ker_{\Z}(M).
\]

Suppose first that $e\in\ker_{\Z}(M)$. Then, on the normalization
$Y$ of the projective closure of each positive-dimensional component,
\[
u^e=u_1^{e_1}\cdots u_k^{e_k}
\]
has zero order along every prime divisor. Thus both $u^e$ and
$(u^e)^{-1}$ have nonnegative order along every prime divisor.
Since $Y$ is normal, both are regular on $Y$; see
\cite[Tag 0AVB]{Stacks}. Hence $u^e$ is a global invertible
regular function on $Y$.
Since $Y$ is proper and integral over $K$, the ring of 
global regular functions $\Gamma(Y,\mathcal{O}_Y)$ is a 
finite-dimensional $K$-algebra; see \cite[Tag 0BUG]{Stacks}. 
Since $Y$ is integral, $\Gamma(Y,\mathcal{O}_Y)$ is an integral domain, 
and a finite-dimensional integral domain over a field is itself a field. 
Hence $\Gamma(Y,\mathcal{O}_Y)$ is a finite field extension of $K$,
and in particular every element of $\Gamma(Y,\mathcal{O}_Y)$ 
is algebraic over $K$. Since $u^e$ is a global regular function on $Y$, 
its image in $K(Y)\simeq\Frac(A/\mathfrak{p})$ is algebraic over $K$.
For a zero-dimensional component, the same conclusion
is automatic, since its function field is a finite extension
of $K$.

For each minimal prime $\mathfrak p$ of $A$, choose a nonzero
polynomial $f_{\mathfrak p}\in K[T]$ annihilating the image of
$u^e$ in $A/\mathfrak p$. With
\[
f=\prod_{\mathfrak p\in\Min(A)}f_{\mathfrak p},
\]
we have
\[
f(u^e)\in
\bigcap_{\mathfrak p\in\Min(A)}\mathfrak p
=\sqrt{0}.
\]
Hence $f(u^e)$ is nilpotent, so
\[
f(u^e)^N=0
\]
for some $N\geq1$. Thus the nonzero polynomial
$f(T)^N\in K[T]$ annihilates $u^e$, and therefore $u^e$ is
algebraic over $K$.
Thus $e\in \Alg_K(u)$.

Conversely, suppose that $e\in \Alg_K(u)$, so there exists a nonzero polynomial $g\in K[T]$ with $g(u^e)=0$ in $A$. Let $\mathfrak{p}$ be any minimal prime of $A$ corresponding to a positive-dimensional component, and let $\pi\colon A\to A/\mathfrak{p}$ be the quotient map. 
Since $u^e$ is a unit in $A$, its image $\pi(u^e)$ is a unit in $A/\mathfrak{p}$, hence nonzero.
Applying $\pi$ to the relation $g(u^e)=0$ gives $g(\pi(u^e))=0$
in $A/\mathfrak{p}$, a nontrivial polynomial relation over $K$
since $g\in K[T]$ is nonzero. 
Hence the image of $u^e$ in $K(Y)\simeq \Frac(A/\mathfrak{p})$ is algebraic over $K$.
Every divisorial valuation $\ord_E$ is trivial on $K$, and hence
also on every element of $K(Y)$ that is algebraic over $K$.
Indeed, choose a polynomial relation
\[
c_0+c_1a+\cdots+c_da^d=0,
\qquad c_0c_d\ne0,
\]
over $K$. If $\ord_E(a)>0$, then the
constant term has uniquely minimal valuation, while if
$\ord_E(a)<0$, then the highest-degree term has uniquely minimal
valuation. Both are impossible.
Therefore
\[
\ord_E(u^e)=0
\]
for every such $E$, so $e\in\ker_{\Z}(M)$. Thus
\[
\Alg_K(u)=\ker_{\Z}(M).
\]

It remains to bound a basis of the integer kernel.
Let $\rho=\rank M$. If $\rho=0$, then
$\ker_{\Z}(M)=\Z^k$, for which the standard basis satisfies the
required bound. If $\rho=k$, then $\ker_{\Z}(M)=0$, which has an empty basis;
the bound holds vacuously. Thus assume $0<\rho<k$.
Choose $\rho$ linearly independent rows of $M$. They span the
row space of $M$ over $\Q$, and hence define the same kernel
over $\Q$, and therefore the same kernel in $\Z^k$.
Choose $\rho$ columns such that the corresponding
$\rho\times\rho$ minor is nonsingular. After permuting the
columns, write the resulting matrix as
\[
(C\;T),
\]
where $C\in\Z^{\rho\times\rho}$ is nonsingular. Put
\[
\Delta=\det C\ne0.
\]
Since all entries have absolute value at most $V$,
\[
|\Delta|\leq \rho!V^\rho.
\]

Write a kernel vector as $(x,y)$, where
$x\in\Z^\rho$ and $y\in\Z^{k-\rho}$. Then
\[
Cx+Ty=0,
\]
and hence
\[
x=-C^{-1}Ty
  =-\frac{\adj(C)Ty}{\Delta}.
\]
Thus projection onto the last $k-\rho$ coordinates identifies
$\ker_{\Z}(M)$ with the full-rank lattice
\[
L_y=
\Set{y\in\Z^{k-\rho}}
{\adj(C)Ty\equiv0\pmod{\Delta}}.
\]
In particular,
\[
\Delta\Z^{k-\rho}\subseteq L_y\subseteq\Z^{k-\rho}.
\]

Let $H_y=(h_{ij})$ be the Hermite normal form of $L_y$, in
the lower-triangular convention
\[
h_{ii}>0,\qquad
h_{ij}=0\quad(i<j),\qquad
0\leq h_{ij}<h_{jj}\quad(i>j);
\]
see \cite[Algorithm~2.4.5]{Cohen93}.
Since
\[
\Delta\Z^{k-\rho}\subseteq L_y,
\]
we have $\Delta e_j\in L_y$ for every $j$.
Expressing $\Delta e_j$ as an integer linear
combination of the rows of $H_y$ and comparing
coordinates from position $k-\rho$ down to position $j$,
triangularity implies
that the coefficients of rows $j+1,\ldots,k-\rho$ vanish.
Consequently, $\Delta$ is an integer multiple of
$h_{jj}$, so
\[
h_{jj}\mid\Delta.
\]
Hence $h_{jj}\leq|\Delta|$, and the reduction condition gives
$|h_{ij}|<h_{jj}\leq|\Delta|$ for every off-diagonal entry.
Thus the rows $y_1,\ldots,y_{k-\rho}$ of $H_y$ form a basis of
$L_y$ satisfying
\[
\|y_j\|_\infty\leq|\Delta|
\leq\rho!V^\rho.
\]

For each $y_j$, let
\[
x_j=-\frac{\adj(C)Ty_j}{\Delta}.
\]
Since
\[
\|Ty_j\|_\infty
\leq(k-\rho)V|\Delta|
\]
and every entry of $\adj(C)$ is a $(\rho-1)\times(\rho-1)$
minor of $C$, of absolute value at most
\[
(\rho-1)!V^{\rho-1},
\]
each coordinate of $\adj(C)Ty_j$ has absolute value at most
\[
\rho(\rho-1)!V^{\rho-1}(k-\rho)V|\Delta|.
\]
Therefore
\[
\|x_j\|_\infty
\leq
\rho!(k-\rho)V^\rho.
\]
Thus
\[
\|(x_j,y_j)\|_\infty
\leq
\rho!(k-\rho)V^\rho
\leq k!V^k
=B_0,
\]
since $\rho!(k-\rho)\leq\rho!(k-\rho)!\leq k!$ and
$V^\rho\leq V^k$ as $V\geq1$.
Since projection onto the last $k-\rho$ coordinates is an
isomorphism from $\ker_{\Z}(M)$ onto $L_y$, the lifted vectors
$(x_j,y_j)$ form a basis of $\ker_{\Z}(M)=\Alg_K(u)$.
This proves the proposition.
\end{proof}

\begin{proposition}\label{pr:alg-lattice2}
A basis of $\Alg_K(u)$ can be computed in deterministic polynomial time.
\end{proposition}

\begin{proof}
Enumerate all vectors
\[
e\in[-B_0,B_0]^k\cap\Z^k.
\]
Since $m$ and $k$ are fixed and $B_0$ is polynomially bounded
in $D$, the number of such vectors is polynomial in the input size.

For each $e=(e_1,\ldots,e_k)$, use the given polynomial
representatives of $u_i$ and $u_i^{-1}$ to construct a polynomial
representative $g_e\in K[X_1,\ldots,X_m]$ of
\[
u^e=u_1^{e_1}\cdots u_k^{e_k}.
\]
Since at most $kB_0$ factors are multiplied,
\[
\deg g_e\leq kB_0D.
\]
As $m$ is fixed, the number of monomials of $g_e$
is polynomially bounded. Its coefficients also have
polynomially bounded bit length over $\Q$;
over $\F_p$, they are reduced modulo $p$.
Thus $g_e$ can be constructed in deterministic
polynomial time.

Apply Lemma~\ref{lem:algebraicity-test} to $g_e+I\in A$ to
determine in deterministic polynomial time whether
$e\in\Alg_K(u)$.
By Proposition~\ref{pr:alg-lattice}, the accepted vectors contain
a basis of $\Alg_K(u)$. Since every accepted vector belongs to
$\Alg_K(u)$, their integer span is exactly $\Alg_K(u)$. A basis of this
integer span can be computed in polynomial time using Hermite
normal form.
\end{proof}

\begin{proof}[Proof of Theorem~\ref{thm:algebraicity-main}]
Proposition~\ref{pr:alg-lattice} gives
the asserted bound on a basis of $\Alg_K(u)$,
while Proposition~\ref{pr:alg-lattice2}
provides a deterministic polynomial-time
algorithm for computing such a basis.
\end{proof}

\subsection{Relation lattices in affine $\Q$-algebras}

We now compute multiplicative relation lattices
in affine $\Q$-algebras.

\begin{proposition}\label{prop:rational-exp}
For fixed $m,k$, the lattice $\Rel_A(u_1,\ldots,u_k)$
can be computed in deterministic polynomial time for specified
units $u_1,\ldots,u_k$ in
\[
A=\Q[X_1,\ldots,X_m]/I,
\qquad
I=(f_1,\ldots,f_s),
\]
with the degree and coefficient encoding above.
\end{proposition}

\begin{proof}
The case $A=0$ is immediate, so assume $A\ne0$. Apply
Proposition~\ref{pr:alg-lattice2} to compute a basis
$b_1,\ldots,b_r$ of $\Alg_\Q(u)$, and put
\[
v_j=u^{b_j}\qquad(1\leq j\leq r).
\]
If $r=0$, then $\Rel_A(u_1,\ldots,u_k)=0$, and we are done.
Otherwise, for each $j$, construct a polynomial representative
of $v_j$ using the given representatives of $u_i$ and $u_i^{-1}$.
Apply Lemma~\ref{lem:algebraicity-test} to compute a monic polynomial
$g_j(T)\in\Q[T]$ annihilating $v_j$, and put
\[
d_j=\deg g_j.
\]
The polynomial representatives of the $v_j$ have polynomially
bounded degree and coefficient bit length. Moreover,
Lemma~\ref{lem:algebraicity-test} gives polynomial bounds on the
degrees and coefficient bit lengths of the $g_j$.

The subalgebra
\[
B=\Q[v_1,\ldots,v_r]\subseteq A
\]
is spanned over $\Q$ by the monomials
\[
w_\alpha=v_1^{\alpha_1}\cdots v_r^{\alpha_r},
\qquad
0\leq\alpha_j<d_j,
\]
since the monic relation $g_j(v_j)=0$ allows every power
$v_j^{n}$ with $n\geq d_j$ to be reduced to a $\Q$-linear
combination of lower powers. Thus this spanning family consists of
\[
M=\prod_{j=1}^r d_j
\]
monomials. Since $r\leq k$ is fixed and each $d_j$ is
polynomially bounded, $M$ is polynomially bounded as well.
Polynomial representatives of all the $w_\alpha$ can be
constructed with polynomially bounded degree and coefficient
bit length.

We next compute the structure constants of $B$. First, we determine
all linear relations among the spanning monomials $w_\alpha$.
For constants $c_\alpha\in\Q$, a relation
\[
\sum_\alpha c_\alpha w_\alpha=0\qquad\text{in }A
\]
is equivalent to an identity
\[
\sum_\alpha c_\alpha W_\alpha
=
\sum_{i=1}^s h_i f_i,
\]
where $W_\alpha\in\Q[X_1,\ldots,X_m]$ is a polynomial
representative of $w_\alpha$ and
$h_i\in\Q[X_1,\ldots,X_m]$.

Let $E$ bound the degrees of the $W_\alpha$.
For any coefficients $c_\alpha\in\Q$, the polynomial
\[
F=\sum_\alpha c_\alpha W_\alpha
\]
has degree at most $E$. Since $m$ is fixed,
Hermann's ideal-membership degree bound
\cite[Theorem~3.4]{Asch} gives a bound polynomial
in $E$ and $D$, independent of the $c_\alpha$,
such that whenever $F\in I$, there exist
polynomials $h_i$ of degrees within this bound
satisfying
\[
F=\sum_{i=1}^s h_i f_i.
\]
As in the proof of Lemma~\ref{lem:bounded-intersection}, we introduce
the coefficients of the $h_i$ as unknowns and compare coefficients
of monomials. This gives a homogeneous linear system over $\Q$
whose projection onto the $c_\alpha$-coordinates is precisely the
space of linear relations among the $w_\alpha$.
Gaussian elimination then selects a basis
\[
\mathcal B=\{w_{\alpha_1},\ldots,w_{\alpha_n}\}
\]
of $B$ from the spanning monomials.

We next compute the multiplication table of $B$ with respect to
the basis $\mathcal B$. For each pair
$w_{\alpha_i},w_{\alpha_k}$ of basis elements, multiply them as
monomials in the $v_j$. Using the relations
\[
g_j(v_j)=0,
\]
reduce all exponents of the $v_j$ to be less than $d_j$. This
expresses the product as a $\Q$-linear combination of the spanning
monomials $w_\alpha$. 
Gaussian elimination also expresses each spanning
monomial $w_\alpha$ as a $\Q$-linear combination
of the basis elements in $\mathcal B$.
Using these expressions, we obtain the structure
constants of $B$.
Their number and coefficient bit lengths are
polynomially bounded. The same linear reduction
gives the coordinates of $1$ and of each $v_j$
in the basis $\mathcal B$.

Each $v_j$ is a unit of $B$. Indeed, multiplication by $v_j$ is
injective on $B$, since $v_j$ is a unit in $A$. As $B$ is
finite-dimensional over $\Q$, this map is also surjective, and hence
$v_j^{-1}\in B$.
We may therefore apply the Lenstra--Silverberg algorithm
\cite[Theorem~1.11(i)]{LS}, which computes in deterministic polynomial
time generators for
\[
\Gamma=
\ker\!\left(
\Z^r\longrightarrow B^\times,\qquad
z\longmapsto\prod_{j=1}^r v_j^{z_j}
\right).
\]
The algorithm applies to arbitrary finite-dimensional commutative
$\Q$-algebras, so no reducedness assumption on $B$ is required.

Let $U\in\Z^{k\times r}$ be the integer matrix whose columns are
$b_1,\ldots,b_r$. Every relation
$e\in\Rel_A(u_1,\ldots,u_k)$ belongs to $\Alg_\Q(u)$, since
$u^e=1$ is algebraic over $\Q$. Hence there is a unique
$z\in\Z^r$ such that
\[
e=Uz.
\]
By the definition of the $v_j$,
\[
u^{Uz}=\prod_{j=1}^r v_j^{z_j},
\]
and, since $B\subseteq A$,
\[
u^{Uz}=1\text{ in }A
\quad\Longleftrightarrow\quad
\prod_{j=1}^r v_j^{z_j}=1\text{ in }B
\quad\Longleftrightarrow\quad
z\in\Gamma.
\]
Therefore
\[
\Rel_A(u_1,\ldots,u_k)=U\Gamma.
\]
Thus, multiplying the computed generators of $\Gamma$ by $U$
gives a generating set for $\Rel_A(u_1,\ldots,u_k)$. A basis
can then be computed in polynomial time using Hermite normal form.
\end{proof}

\section{Relation lattices over rings of prime-power characteristic}
\label{sec:primepower}

In this section we prove that a basis of the relation lattice
$\Rel_A(u_1,\ldots,u_k)$
can be computed in bounded-error quantum polynomial time
for specified units $u_1,\ldots,u_k$ in a ring of the form
\[
A=(\Z/p^a\Z)[X_1,\ldots,X_m]/I,
\qquad
I=(f_1,\ldots,f_s),
\]
where $p$ is prime and $a\geq 1$.
We assume that the number $m$ of polynomial variables and the
number $k$ of specified units are fixed.
The input consists of the generators $f_1,\ldots,f_s$ of $I$,
polynomial representatives of the units $u_1,\ldots,u_k$
and their inverses, and the modulus $p^a$ in binary.
In particular, the input size includes the bit length
of $p^a$, which is $\lfloor a\log_2 p\rfloor+1$.

We first show how to construct an exact presentation of
a finite subring generated by elements satisfying monic
polynomial relations.

Let $0\le r\le k$, and suppose that
$v_1,\ldots,v_r\in A$ are units given by polynomial
representatives of polynomial size.
Suppose further that, for each $j$,
a monic polynomial $g_j(T)\in(\Z/p^a\Z)[T]$
of polynomially bounded degree is given such that
$g_j(v_j)=0$.
Put
\[
B=(\Z/p^a\Z)[v_1,\ldots,v_r]\subseteq A.
\]

\begin{lemma}[An exact finite-subring presentation]
\label{lem:finite-subring}
The ring $B$ has an exact additive presentation and
multiplication table computable in deterministic polynomial
time. In particular, $B$ admits unique polynomial-length
encodings and polynomial-time ring operations.
\end{lemma}

\begin{proof}
Put $d_j=\deg g_j$. As in the proof of
Proposition~\ref{prop:rational-exp}, the monic relations
$g_j(v_j)=0$ imply that $B$ is generated as a
$\Z/p^a\Z$-module by the monomials
\[
w_\alpha=v_1^{\alpha_1}\cdots v_r^{\alpha_r},
\qquad
0\leq\alpha_j<d_j.
\]
The number of these monomials is
\[
M=\prod_{j=1}^r d_j.
\]
Since $r\le k$ and $k$ is fixed, the polynomial bounds
on the $d_j$ imply that $M$ is polynomially bounded.
Polynomial representatives
of all the $w_\alpha$ can be constructed with polynomially bounded
degree and coefficient bit length.

We next determine all additive relations among the spanning
monomials. Unlike the rational case, these relations form an
integer lattice rather than a vector space over $\Q$.
Enumerate the spanning monomials as $w_1,\ldots,w_M$ and put
\[
\Omega=
\Set{c\in\Z^M}{\sum_{i=1}^M c_iw_i=0\text{ in }A}.
\]
Then
\[
B_{\rm add}\simeq\Z^M/\Omega.
\]

For each $i$, choose an integer polynomial
$W_i\in\Z[X_1,\ldots,X_m]$ lifting a polynomial representative
of $w_i$, and choose integer lifts, still denoted by
$f_1,\ldots,f_s$, of the generators of $I$.
Thus
\[
A\simeq
\Z[X_1,\ldots,X_m]/(f_1,\ldots,f_s,p^a).
\]

Every $c\in \Omega$ is congruent modulo $p^a\Z^M$ to a vector
$c^{(0)}\in \Omega$ satisfying
\[
0\leq c_i^{(0)}<p^a.
\]
For such a vector, the polynomial
\[
\sum_{i=1}^M c_i^{(0)}W_i
\]
has polynomially bounded degree and coefficient bit length.

By Aschenbrenner's integer ideal-membership bound
\cite[Theorem~A]{Asch}, there is a uniform bound
$E$, polynomial in the input size for fixed $m$,
such that every identity
\[
\sum_{i=1}^M c_i^{(0)}W_i
=
\sum_{\nu=1}^s f_\nu h_\nu+p^ah_0,
\qquad 0\leq c_i^{(0)}<p^a,
\]
that holds for some integer polynomials $h_\nu$
also holds with $\deg h_\nu\leq E$ for all $\nu$.
Here the bound is uniform because the degrees
and logarithmic heights of the polynomials on
the left are polynomially bounded.

We now consider the equation
\begin{equation}\label{eq:finite-relation}
\sum_{i=1}^M c_iW_i
=
\sum_{\nu=1}^s f_\nu h_\nu+p^ah_0,
\end{equation}
where the $c_i$ are unknown integer constants and
$h_0,h_1,\ldots,h_s$ are unknown integer polynomials
of degree at most $E$.
Expanding the $h_\nu$ and comparing coefficients of monomials
gives a homogeneous integer linear system with polynomially
many equations and unknowns and polynomially bounded
coefficient bit lengths.

The projection of its integer solution lattice onto the
$c$-coordinates is exactly $\Omega$. Indeed, every projected
solution belongs to $\Omega$. Conversely, the ideal-membership
degree bound shows that every relation
$c^{(0)}\in\Omega$ with $0\le c_i^{(0)}<p^a$
arises from a solution of this system. The same holds for
$p^ae_i$, by taking $h_0=W_i$ and all other $h_\nu=0$.
Every element of $\Omega$ is the sum of a bounded relation
$c^{(0)}$ and an element of $p^a\Z^M$.
Hence the projection is all of $\Omega$.

Standard polynomial-time integer linear algebra
computes a basis of the integer solution lattice.
Projecting this basis onto the $c$-coordinates
and applying Hermite normal form gives a basis
of $\Omega$ in deterministic polynomial time.

Smith normal form, together with the corresponding 
unimodular transformation matrices, now gives an explicit decomposition
\[
B_{\rm add}=\Z^M/\Omega
\simeq
\Z/n_1\Z\oplus\cdots\oplus\Z/n_t\Z.
\]
Since $p^a\Z^M\subseteq \Omega$, each $n_i$ divides $p^a$.
The standard representatives $0,\ldots,n_i-1$ in each
cyclic factor therefore give unique encodings of polynomial
length.

We next compute the multiplication table of $B$.
For each pair $w_\alpha,w_\beta$, their product is a
monomial in the $v_j$. Using the monic relations
$g_j(v_j)=0$, we express this product as a
$\Z/p^a\Z$-linear combination of the spanning
monomials. The unimodular transformations obtained
from Smith normal form then allow us to express
products of the canonical additive generators
in Smith coordinates.

The resulting multiplication table has polynomial
size, and ring operations in the canonical
representation can be performed in deterministic
polynomial time.
\end{proof}

Given a group $G$, a subgroup $H\le G$, and a set 
$X$, we say a function $f:G\to X$ \emph{hides} the subgroup $H$
if for all 
$g_{1},g_{2}\in G,f(g_{1})=f(g_{2})$
if and only if 
$g_{1}H=g_{2}H$.
Equivalently, 
$f$
is constant on each coset of $H$, 
while it is different between the different cosets of $H$.

\begin{lemma}[Relations in a finite ring]
\label{lem:finitekernel}
There is a bounded-error quantum polynomial-time algorithm
that computes a basis of $\Rel_B(v_1,\ldots,v_r)$.
\end{lemma}

\begin{proof}
By Lemma~\ref{lem:finite-subring}, the ring $B$
has unique polynomial-length encodings and
efficient ring operations. Since the $v_i$ are units,
we can use Shor's order-finding algorithm \cite{Shor}
to compute
\[
o_i=\ord_{B^\times}(v_i)\le |B|,
\qquad i=1,\ldots,r.
\]
Since $|B|\le(p^a)^M$, the orders $o_i$ have
polynomially bounded bit lengths.

Consider the composition of group homomorphisms
\[
\Z^r
\stackrel{\varphi}{\longrightarrow}
\bigoplus_{i=1}^r\Z/o_i\Z
\stackrel{\psi}{\longrightarrow}B^\times,
\]
where $\varphi$ is the natural projection and
$\psi(e)=v^e$. By definition,
\[
\Rel_B(v_1,\ldots,v_r)=\ker(\psi\circ\varphi)
=\varphi^{-1}(\ker\psi).
\]

The homomorphism $\psi$ hides the subgroup $\ker\psi$.
The quantum algorithm for the finite abelian hidden subgroup
problem \cite{CM} therefore computes generators
$g_1,\ldots,g_t$ of $\ker\psi$ in bounded-error
polynomial time.

Choose lifts $\widetilde g_1,\ldots,\widetilde g_t\in\Z^r$
of these generators. Then
\[
\Rel_B(v_1,\ldots,v_r)
=
\left\langle
\widetilde g_1,\ldots,\widetilde g_t,
o_1e_1,\ldots,o_re_r
\right\rangle,
\]
where $e_1,\ldots,e_r$ are the standard basis vectors
of $\Z^r$.
Applying Hermite normal form gives a basis of
$\Rel_B(v_1,\ldots,v_r)$ in deterministic polynomial time.
\end{proof}

\begin{proposition}\label{prop:primepower}
Under the standing assumptions of this section,
a basis of $\Rel_A(u_1,\ldots,u_k)$
can be computed in bounded-error quantum polynomial time.
\end{proposition}

\begin{proof}
The case $A=0$ is immediate, so assume $A\ne0$.
Put $\overline A=A/pA$, and let $\overline u_i$ denote
the image of $u_i$ in $\overline A$.
Apply Theorem~\ref{thm:algebraicity-main} to compute a basis
$b_1,\ldots,b_r$ of
\[
L=\Alg_{\F_p}(\overline u_1,\ldots,\overline u_k).
\]
If $r=0$, then $\Rel_A(u_1,\ldots,u_k)=0$, and we are done.

Put $v_j=u^{b_j}$ for $j=1,\ldots,r$.
By Proposition~\ref{pr:alg-lattice2}, the coordinates
of the $b_j$ are polynomially bounded.
Using the given polynomial representatives of the $u_i$
and their inverses, we can therefore construct polynomial
representatives of the $v_j$ with polynomially bounded
degree and coefficient bit length.
Apply Lemma~\ref{lem:algebraicity-test} to each
$\overline v_j\in\overline A$ to compute a monic polynomial
$g_j(T)\in\F_p[T]$ satisfying
\[
g_j(\overline v_j)=0.
\]
The degrees and coefficient bit lengths of these polynomials
are polynomially bounded.

Choose a monic coefficient lift
$\widetilde g_j(T)\in(\Z/p^a\Z)[T]$ of $g_j(T)$.
Since $g_j(\overline v_j)=0$, we have
\[
\widetilde g_j(v_j)\in pA.
\]
Consequently,
\[
\widetilde g_j(v_j)^a\in p^aA=0.
\]
Thus
\[
G_j(T)=\widetilde g_j(T)^a
\]
is a monic polynomial annihilating $v_j$, with $\deg G_j=a\deg g_j$.
Since $a\leq\log_2(p^a)$, the polynomials $G_j$
have polynomially bounded degrees and can be
constructed in deterministic polynomial time.

By Lemma~\ref{lem:finite-subring}, the subring
\[
B=(\Z/p^a\Z)[v_1,\ldots,v_r]\subseteq A
\]
has an exact finite presentation with unique polynomial-length
encodings and polynomial-time ring operations.

Each $v_j$ is a unit of $B$. Indeed, multiplication by $v_j$
is injective on $B$, since $v_j$ is a unit in $A$. As $B$ is
finite, this map is also surjective, and hence $v_j^{-1}\in B$.

We may therefore apply Lemma~\ref{lem:finitekernel} to compute
a basis of
\[
\Gamma=\Rel_B(v_1,\ldots,v_r)
\]
in bounded-error quantum polynomial time.

Let $U\in\Z^{k\times r}$ be the matrix with columns
$b_1,\ldots,b_r$.
Since
\[
\Rel_A(u_1,\ldots,u_k)
\subseteq
\Rel_{\overline A}(\overline u_1,\ldots,\overline u_k)
\subseteq L,
\]
every relation $e\in\Rel_A(u_1,\ldots,u_k)$
has a unique representation $e=Uz$ with $z\in\Z^r$.
By the definition of the $v_j$,
\[
u^{Uz}=\prod_{j=1}^r v_j^{z_j},
\]
and, since $B\subseteq A$,
\[
u^{Uz}=1\text{ in }A
\quad\Longleftrightarrow\quad
\prod_{j=1}^r v_j^{z_j}=1\text{ in }B
\quad\Longleftrightarrow\quad
z\in\Gamma.
\]
Therefore
\begin{equation}\label{eq:primepower-transport}
\Rel_A(u_1,\ldots,u_k)=U\Gamma.
\end{equation}
Multiplying the computed basis vectors of $\Gamma$ by $U$
gives a basis of $\Rel_A(u_1,\ldots,u_k)$.
\end{proof}

\section{Proof of Theorem~\ref{thm:affine-main}}
\label{sec:integral-exp}

Let $S=\Z[X_1,\ldots,X_m]/I$.
We first determine whether $S$ is the zero ring
by testing whether $1\in I$.
By Proposition~\ref{prop:linear}, this can be done
in deterministic polynomial time.
If $S=0$, then every product
$u_1^{e_1}\cdots u_k^{e_k}$ equals $1$ in $S$
for every $e=(e_1,\ldots,e_k)\in\Z^k$.
Consequently,
\[
\Rel_S(u_1,\ldots,u_k)=\Z^k,
\]
and we output the standard basis of $\Z^k$.
Henceforth assume that $S\ne0$.

Compute the integer $c$ from
Proposition~\ref{th:arithmetic} and factor it
using Shor's quantum integer-factorization
algorithm \cite{Shor}. Write
\[
c=\prod_{\nu=1}^t p_\nu^{a_\nu}.
\]
Since $c$ has polynomially bounded bit length,
the number $t$ of prime factors and the bit lengths
of the prime powers $p_\nu^{a_\nu}$ are polynomially
bounded.
By Proposition~\ref{th:arithmetic}, the natural map
\[
S\longrightarrow (S\otimes\Q)\times S/cS
\]
is injective. Moreover, the Chinese remainder theorem gives
\[
S/cS\simeq
\prod_{\nu=1}^t S/p_\nu^{a_\nu}S.
\]
Consequently,
\begin{equation}\label{eq:integral-intersection}
\Rel_S(u_1,\ldots,u_k)
=
\Rel_{S\otimes\Q}(u_1,\ldots,u_k)
\cap
\bigcap_{\nu=1}^t
\Rel_{S/p_\nu^{a_\nu}S}(u_1,\ldots,u_k).
\end{equation}

The relation lattice over $S\otimes\Q$ is computed by
Proposition~\ref{prop:rational-exp}, and those over the
prime-power quotients are computed by
Proposition~\ref{prop:primepower}.
Zero-ring factors impose no restriction.
If $c=1$, the intersection over the prime-power
quotients is omitted.

It remains to compute the intersection of the resulting
integer lattices. Let $B_1,\ldots,B_\ell$ be their basis
matrices, where $B_i$ has $r_i$ linearly independent
columns. Introduce a common vector $x\in\Z^k$ and
require
\[
x=B_i z_i,\qquad z_i\in\Z^{r_i},
\quad i=1,\ldots,\ell.
\]
These equations form a homogeneous integer linear system
whose solution lattice projects onto the required
intersection.

Since each $B_i$ has linearly independent columns,
the projection onto the $x$-coordinates is injective.
Thus, projecting a basis of the integer solution
lattice gives a basis of the intersection.

The combined system has polynomially many rows and
columns, and its coefficients have polynomially bounded
bit lengths. Standard polynomial-time integer kernel
algorithms therefore compute a basis of its solution
lattice with polynomially bounded bit lengths.
The same bounds hold for the projected basis.

All remaining classical computations take
deterministic polynomial time.
Amplifying the polynomially many quantum subroutines
ensures that the overall probability of correctness
is at least $2/3$.
Thus a basis of $\Rel_S(u_1,\ldots,u_k)$ can be
computed in bounded-error quantum polynomial time.
\hfill$\square$

\section{Proof of Corollary~\ref{cor:affine-membership}}
\label{sec:affine-membership}

Apply Theorem~\ref{thm:affine-main} to the $k+1$ units
$h,u_1,\ldots,u_k$ to compute a basis
$b_1,\ldots,b_r$ of
\[
\Lambda=\Rel_S(h,u_1,\ldots,u_k)
\leq\Z^{k+1}.
\]
Write
\[
b_j=(c_j,v_j),
\qquad c_j\in\Z,\quad v_j\in\Z^k,
\]
and put
\[
g=\gcd(c_1,\ldots,c_r),
\]
where $g=0$ if $r=0$ or all $c_j$ are zero.

The subgroup of first coordinates of elements of
$\Lambda$ is $g\Z$.
Since
\[
(-1,e)\in\Lambda
\quad\Longleftrightarrow\quad
h^{-1}u^e=1
\quad\Longleftrightarrow\quad
u^e=h,
\]
we have
\[
h\in\langle u_1,\ldots,u_k\rangle
\quad\Longleftrightarrow\quad
-1\in g\Z
\quad\Longleftrightarrow\quad
g=1.
\]

Suppose that $g=1$. By the extended Euclidean
algorithm, we can find integers
$\lambda_1,\ldots,\lambda_r$ such that
\[
\sum_{j=1}^r\lambda_jc_j=-1.
\]
Put
\[
e_0=\sum_{j=1}^r\lambda_jv_j.
\]
Then $(-1,e_0)\in\Lambda$, so $u^{e_0}=h$.
Moreover, for any $e\in\Z^k$,
\[
u^e=h
\quad\Longleftrightarrow\quad
u^{e-e_0}=1
\quad\Longleftrightarrow\quad
e-e_0\in\Rel_S(u_1,\ldots,u_k).
\]
Consequently,
\[
\{e\in\Z^k:u^e=h\}
=e_0+\Rel_S(u_1,\ldots,u_k).
\]
All the additional computations use polynomial-time
integer linear algebra. Since $k+1$ is fixed,
the application of Theorem~\ref{thm:affine-main}
also takes bounded-error quantum polynomial time.
\hfill$\square$

\section{The conjugacy-to-orbit reduction}
\label{sec:group}

In this section, we reduce the conjugacy problem in a
finitely generated metabelian group to an orbit problem
in a finitely presented module over its integral group ring.
We also show that this reduction can be carried out
in polynomial time.

Let $G$ be a $k$-generated metabelian group given by
\eqref{eq:G},
where the defining relators are part of the input. Let
\[
A=G',\qquad Q=G/A
\]
be the derived subgroup and the abelianization of $G$,
respectively, and let $\varphi:G\to Q$ be the natural projection.
Since $G$ is metabelian, both $A$ and $Q$ are abelian.

Conjugation by elements of $G$ induces an action of $Q$ on $A$.
Indeed, every element of $A$ acts trivially on $A$ by conjugation.
Thus $A$ is naturally a module over the commutative group ring
\[
R=\mathbb Z Q.
\]
We use right conjugation, $a^g=g^{-1}ag$, and write $A$
additively when regarded as an $R$-module. In particular,
\[
a\cdot\varphi(g)=a^g,\qquad
[a,g]=a\cdot(\varphi(g)-1)
\]
for all $a\in A$ and $g\in G$.

\subsection{$R$-module coordinates in the derived subgroup}
\label{se:coordinates}

The derived subgroup $A=G'$ is generated as an
$R$-module by the commutators
\[
a_{ij}=[x_i,x_j],\qquad 1\leq i<j\leq k.
\]
Thus every $a\in A$ can be written as
\[
a=\sum_{1\leq i<j\leq k} a_{ij}f_{ij},
\qquad f_{ij}\in R.
\]
We call $(f_{ij})_{i<j}$ an
\emph{$R$-module coordinate vector} of $a$.
Such a vector need not be unique because the
generators $a_{ij}$ may satisfy nontrivial
$R$-module relations.

\begin{proposition}\label{pr:collection}
Fix $k\geq 1$.
There is a deterministic algorithm that,
given a $k$-generated metabelian group $G$
defined by \eqref{eq:G} and a word
$w=w(x_1,\ldots,x_k)$ representing an element
of $G'$, produces an $R$-module coordinate
vector $(f_{ij})_{i<j}$ of $w$ satisfying
the following:
\begin{itemize}
\item
The running time is polynomial in $N+|w|$,
where $N=m+\sum_{i=1}^m |r_i|$, with degree
depending only on $k$.
\item
Each coordinate $f_{ij}$ has numerical degree,
number of nonzero terms, and coefficient
bit length bounded by a polynomial in
$N+|w|$, with degree depending only on $k$.
\end{itemize}
\end{proposition}

\begin{proof}
Let $F=F(x_1,\ldots,x_k)$, and let
$\alpha:F\to\Z^k$ be the abelianization map.
Put
\[
v_i=\alpha(r_i)\in\Z^k,
\qquad 1\leq i\leq m,
\]
and let $V$ be the integer matrix with columns
$v_1,\ldots,v_m$.

Since $w$ represents an element of $G'$, we have
\[
\alpha(w)\in\langle v_1,\ldots,v_m\rangle.
\]
Thus the integer system
\[
Vc=-\alpha(w)
\]
has a solution.

The matrix $V$ has $k$ rows, and its entries
have absolute value at most $N$, while
$\|\alpha(w)\|_\infty\leq |w|$.
Since $k$ is fixed, standard bounds for
solvable integer linear systems imply that
there exists a solution $c\in\Z^m$ with
\[
Vc=-\alpha(w),
\qquad
\sum_{i=1}^m|c_i|
\leq (N+|w|)^{O(k)}.
\]
Moreover, such a solution can be computed
in deterministic polynomial time using
integer linear algebra algorithms.

Set
\[
u=r_1^{c_1}\cdots r_m^{c_m}.
\]
The word $u$ represents the identity in $G$,
and its length satisfies
\[
|u|\leq (N+|w|)^{O(k)}.
\]
Moreover,
\[
\alpha(wu)=0,
\]
so $wu\in F'$. Thus $w$ and $wu$ represent
the same element of $G'$.
Put
\[
L=|wu|\leq (N+|w|)^{O(k)}.
\]
Collect the word $wu$ by sorting its letters
into the order
\[
x_1^{\pm1},x_2^{\pm1},\ldots,x_k^{\pm1}.
\]
Whenever adjacent letters
$a=x_j^{\pm1}$ and $b=x_i^{\pm1}$,
with $j>i$, are interchanged, use
\[
sabv=sbav\cdot[a,b]^v.
\]
Each interchange introduces a conjugate
of a commutator of two generator letters.

At most $O(L^2)$ interchanges are required,
and every conjugating word has length $O(L)$.
Each commutator
$[x_j^{\pm1},x_i^{\pm1}]$
is a conjugate of
$[x_i,x_j]^{\pm1}$
by a word of bounded length.

Since $wu\in F'$, the exponent sum of each
generator in the sorted word is zero.
Consequently, the sorted word freely reduces
to the identity.
Hence
\[
wu=\prod_{\ell=1}^{T}
\bigl([x_{i_\ell},x_{j_\ell}]
^{\varepsilon_\ell}\bigr)^{h_\ell},
\]
where $T=O(L^2)$,
$\varepsilon_\ell\in\{1,-1\}$,
and $|h_\ell|=O(L)$.

Passing to $G$ and writing $G'$ additively,
we obtain
\[
w=\sum_{\ell=1}^{T}
\varepsilon_\ell a_{i_\ell j_\ell}
\varphi(h_\ell).
\]
Grouping terms corresponding to each pair
$i<j$ gives an $R$-module coordinate vector
\[
w=\sum_{1\leq i<j\leq k}a_{ij}f_{ij},
\quad\text{where }
f_{ij}=
\sum_{\ell:\,(i_\ell,j_\ell)=(i,j)}
\varepsilon_\ell\varphi(h_\ell).
\]
We represent these coordinates as Laurent
polynomials in the images of
$x_1,\ldots,x_k$ in $Q$, without requiring
canonical representatives of elements of $Q$.

Each $f_{ij}$ has at most $O(L^2)$ nonzero
terms. The absolute values of the exponents of
$x_1,\ldots,x_k$ in these Laurent monomials
are $O(L)$, and every coefficient has
absolute value at most $O(L^2)$.
Since $L\leq (N+|w|)^{O(k)}$, the numerical
degrees and support sizes are bounded by
$(N+|w|)^{O(k)}$, while the coefficient
bit lengths are $O(k\log(N+|w|))$.
All constructions can be performed in
deterministic polynomial time in $N+|w|$,
with degree depending only on $k$.
\end{proof}

\subsection{From conjugacy to orbit membership}
\label{se:noskov}

Let $x,y\in G$ be the input elements of the conjugacy
problem. Since $Q$ is abelian, conjugate elements have
the same image under $\varphi$. We may therefore assume that
\[
\varphi(x)=\varphi(y)=q\in Q.
\]
For every $a\in A$, we have
\[
[a,x]=a\cdot(q-1).
\]
Hence, the subgroup
\[
D=(q-1)A
\]
consists precisely of the commutators $[a,x]$ with $a\in A$.

\begin{lemma}\label{le:D-normal}
The following hold:
\begin{itemize}
\item[(a)] $D\unlhd G$;
\item[(b)] $xD$ centralizes $A/D$ in $G/D$.
\end{itemize}
\end{lemma}

\begin{proof}
Since $D$ is an $R$-submodule of $A$, it is normal
in $G$. Moreover, for every $a\in A$,
\[
[a,x]=a\cdot(q-1)\in D.
\]
Hence $xD$ commutes with every element of $A/D$.
\end{proof}

Put
\[
M=A/D,\qquad C=\langle M,xD\rangle\leq G/D.
\]

\begin{lemma}\label{le:C-module}
The following hold:
\begin{itemize}
\item[(a)] $C$ is an abelian normal subgroup of $G/D$;
\item[(b)] conjugation induces an $R$-module structure on $C$;
\item[(c)] $yD\in C$.
\end{itemize}
\end{lemma}

\begin{proof}
By Lemma~\ref{le:D-normal}(b), the element $xD$
centralizes $M$. Since $M$ is abelian, so is $C$.
For every $g\in G$, we have
\[
x^g=x[x,g],
\]
where $[x,g]\in A$. Thus conjugation by $gD$
preserves $C$, and $C$ is normal in $G/D$.

Since $M\subseteq C$ and $C$ is abelian, conjugation
by elements of $A$ acts trivially on $C$. Therefore
the conjugation action of $G$ on $C$ factors through
$Q=G/A$, making $C$ an $R$-module.

Finally, since $\varphi(x)=\varphi(y)$,
we have $x^{-1}y\in A$. Therefore
\[
yD=(xD)(x^{-1}yD)\in C,
\]
proving (c).
\end{proof}

The following proposition shows that conjugacy in $G$
is equivalent to an orbit problem for the action of $Q$
on $C$.

\begin{proposition}[Noskov reduction]\label{prop:noskov}
Let $x,y\in G$ have the same image in $Q$, and let
$D$ and $C$ be as defined above. Then
\begin{equation}\label{eq:noskov}
x\sim_G y
\quad\Longleftrightarrow\quad
yD=(xD)\cdot s
\quad\text{for some }s\in Q.
\end{equation}
\end{proposition}

\begin{proof}
($\Rightarrow$)
If $x^g=y$ for some $g\in G$, then
\[
yD=(xD)\cdot\varphi(g),
\]
proving the forward implication.

($\Leftarrow$)
Conversely, suppose that
\[
yD=(xD)\cdot s
\]
for some $s\in Q$. Choose $g\in G$ with
$\varphi(g)=s$. Then $x^gD=yD$, so
\[
x^g=yd
\]
for some $d\in D$.

Since $\varphi(y)=q$, the map
$A\longrightarrow A$ defined by
$a\longmapsto[y,a]$
has image $(1-q)A=D$. We may therefore choose
$a\in A$ such that $[y,a]=d^{-1}$. As $A$ is
abelian, $a$ centralizes $d$, and hence
\[
x^{ga}=(x^g)^a=(yd)^a
       =y[y,a]d=y.
\]
Thus $x$ and $y$ are conjugate in $G$.
\end{proof}

\subsection{An explicit $R$-module presentation of $A$}
\label{se:A-presentation}

We now show that a finite $R$-module presentation
of $A=G'$ can be constructed efficiently from
the defining relators of $G$.

\begin{proposition}\label{prop:A-presentation}
Fix $k\geq 1$.
There is a deterministic algorithm that,
given a $k$-generated metabelian group $G$
defined by \eqref{eq:G},
constructs a finite $R$-module presentation
of $A=G'$ with generators
$a_{ij}=[x_i,x_j]$ for $1\leq i<j\leq k$,
in time polynomial in
\[
N=m+\sum_{i=1}^m |r_i|,
\]
with degree depending only on $k$.
The presentation has polynomially many
relations, whose coefficients admit
Laurent polynomial representatives in
$\Z[t_1^{\pm1},\ldots,t_k^{\pm1}]$
with polynomially bounded degrees,
support sizes, and coefficient bit lengths.
\end{proposition}

\begin{proof}

Let $F=F(x_1,\ldots,x_k)$ be the free group,
and put
\[
H=F/F'',\qquad
T=\Z[t_1^{\pm1},\ldots,t_k^{\pm1}].
\]
Then $H$ is the free metabelian group
of rank $k$, and $H'$ is a $T$-module
under conjugation.
The module $H'$ is generated by $a_{ij}=[x_i,x_j]$,
with defining relations
\begin{equation}\label{eq:free-metabelian-relations}
t_it_j(t_\ell-1)a_{ij}
-t_it_\ell(t_j-1)a_{i\ell}
+t_jt_\ell(t_i-1)a_{j\ell}=0,
\qquad 1\leq i<j<\ell\leq k.
\end{equation}
Indeed, after introducing the
normalized generators
\[
b_{ij}=t_it_j a_{ij},
\]
these become the Koszul relations
\[
b_{ij}(t_\ell-1)
-b_{i\ell}(t_j-1)
+b_{j\ell}(t_i-1)=0.
\]
Their completeness follows from
the Magnus embedding and the
exactness of the Koszul complex
associated with
$t_1-1,\ldots,t_k-1$.

Let $K$ be the normal closure of
$r_1,\ldots,r_m$ in $H$. Then
\[
G=H/K,\qquad
A=G'\cong H'/(K\cap H').
\]
Let $\alpha:H\to\Z^k$ be abelianization,
and set
\[
v_i=\alpha(r_i),\qquad
V=(v_1\ \cdots\ v_m)\in\Z^{k\times m}.
\]
Thus
\[
Q\cong\Z^k/\operatorname{im}V,
\qquad
R=\Z Q\cong T/J,
\]
where
\[
J=(t^{v_1}-1,\ldots,t^{v_m}-1),
\qquad
t^v=t_1^{v_1}\cdots t_k^{v_k}.
\]

We first determine generators of
$K\cap H'$ as a $T$-module.
Set
\[
\mathcal S_1=
\{[r_i,x_j]\mid
1\leq i\leq m,\ 1\leq j\leq k\}.
\]
Choose a $\Z$-basis
$c^{(1)},\ldots,c^{(s)}$ of
$\ker_{\Z}V$, and define
\[
\rho_\ell=
r_1^{c_1^{(\ell)}}\cdots
r_m^{c_m^{(\ell)}},
\qquad 1\leq\ell\leq s.
\]
Since $\alpha(\rho_\ell)=Vc^{(\ell)}=0$,
each $\rho_\ell$ belongs to $K\cap H'$.
Put
\[
\mathcal S_2=\{\rho_1,\ldots,\rho_s\}.
\]

We claim that
\begin{equation}\label{eq:K-intersection}
K\cap H'=
\langle\mathcal S_1\cup\mathcal S_2\rangle_T.
\end{equation}
The normal closure of $\mathcal S_1$
in $H$ is $[K,H]$: modulo this normal
closure, every $r_i$ commutes with all
generators $x_j$, and hence $K$ is central.
Since $[K,H]\subseteq H'$ and $H'$ is
abelian, it follows that
\[
[K,H]=\langle\mathcal S_1\rangle_T.
\]

The quotient $K/[K,H]$ is abelian and
generated by the images of the $r_i$.
Consequently, the homomorphism
\[
\beta:\Z^m\longrightarrow K/[K,H],
\qquad
c\longmapsto
r_1^{c_1}\cdots r_m^{c_m}[K,H]
\]
is surjective, and its composition
with the map induced by $\alpha$ is $V$.

For any $z\in K\cap H'$, choose
$c\in\Z^m$ with $z[K,H]=\beta(c)$.
Then $Vc=\alpha(z)=0$, so
\[
c=\sum_{\ell=1}^s n_\ell c^{(\ell)}
\]
for some integers $n_\ell$. Hence
\[
z[K,H]=
\prod_{\ell=1}^s
\rho_\ell^{n_\ell}[K,H].
\]
Since $[K,H]$ is generated by
$\mathcal S_1$ as a $T$-module,
this proves \eqref{eq:K-intersection}.

The conjugation action of $H$ on $A$
factors through $Q=G/A$.
Consequently, the ideal $J$ annihilates
$A$, and hence
\[
JH'\subseteq K\cap H'.
\]
By \eqref{eq:K-intersection}, an
$R$-module presentation of $A$
is obtained from
\eqref{eq:free-metabelian-relations}
by adjoining the coordinate vectors
of the elements of
$\mathcal S_1\cup\mathcal S_2$,
computed over $T$ and mapped
componentwise to $R$.

We now establish the complexity bounds.
The matrix $V$ has $k$ rows, and its
entries have absolute value at most $N$.
Put $r=\operatorname{rank}V\leq k$.
If $r=0$, we take the standard basis
of $\ker_{\Z}V=\Z^m$.
Suppose that $r>0$. Every $r\times r$
minor of $V$ has absolute value at most
\[
D=r!N^r\leq k!N^k.
\]
Choose $r$ independent columns and
$r$ rows forming a nonsingular
$r\times r$ submatrix $B$.
After permuting columns, write the
corresponding row submatrix of $V$
as $(B\ C)$.
Since these rows span the row space
of $V$ over $\Q$, we have
\[
\ker_{\Z}V=
\{(u,w)\in\Z^r\oplus\Z^{m-r}
\mid Bu+Cw=0\}.
\]

Put $d=|\det B|\leq D$.
The possible vectors $w$ form the
sublattice
\[
L=\{w\in\Z^{m-r}\mid
\operatorname{adj}(B)Cw
\equiv 0\pmod d\}.
\]
This lattice is the projection onto
the first $m-r$ coordinates of the
integer solution space of
\[
\operatorname{adj}(B)Cw-dz=0,
\qquad
(w,z)\in\Z^{m-r}\oplus\Z^r.
\]
Since $z$ is uniquely determined by
$w$, this projection is an isomorphism.
Thus $L$ can be computed in deterministic
polynomial time using integer normal
form algorithms.

Moreover, $L$ is the kernel of a
homomorphism
\[
\Z^{m-r}\longrightarrow(\Z/d\Z)^r,
\]
so
\[
[\Z^{m-r}:L]\leq d^r\leq D^k.
\]
A reduced Hermite normal form basis
of $L$ therefore has entries bounded
by $D^k$.
For each basis vector $w$ of $L$, put
\[
u=-B^{-1}Cw.
\]
The defining congruence of $L$
ensures that $u$ is integral.
By Cramer's rule, each entry of
$B^{-1}C$ is a ratio of two
$r\times r$ minors of $V$, and hence
has absolute value at most $D$.
Consequently,
\[
\|(u,w)\|_1
\leq (m-r)D^k+r(m-r)D^{k+1}
\leq N^{O(k^2)}.
\]
These vectors form a $\Z$-basis
$c^{(1)},\ldots,c^{(s)}$ of
$\ker_{\Z}V$, where $s=m-r\leq m$.
It follows that
\[
|\rho_\ell|
\leq
\sum_{i=1}^m
|c_i^{(\ell)}|\,|r_i|
\leq N^{O(k^2)}.
\]

The words $[r_i,x_j]$ have length
$O(|r_i|+1)$, and the words $\rho_\ell$
have polynomial length.
Since they all belong to $H'$, their
$T$-module coordinate vectors can be
computed in deterministic polynomial
time by Proposition~\ref{pr:collection}.
Their degrees, support sizes, and
coefficient bit lengths are polynomially
bounded in $N$, with degree depending
only on $k$.

Finally, the presentation contains
at most
\[
\binom{k}{3}+mk+s
\leq \binom{k}{3}+mk+m
\]
relations. Hence it has polynomial
size and can be constructed in
deterministic polynomial time,
with degree depending only on $k$.
\end{proof}

\subsection{An explicit $R$-module presentation of $C$}
\label{se:C-presentation}

Noskov's reduction is the structural step underlying his
decidability theorem \cite{Noskov}. To obtain a polynomial-time
reduction, we must construct an explicit $R$-module presentation
of the input-dependent module $C$.

By Proposition~\ref{prop:A-presentation}, a finite
$R$-module presentation
\[
A=\left\langle
a_{ij}\ (1\leq i<j\leq k)
\ \middle|\ \mathcal R_A
\right\rangle_R
\]
can be constructed in polynomial time from the
defining relations of $G$.

Recall that $M=A/(q-1)A$. Its presentation
is obtained by adjoining the relations
\[
(q-1)a_{ij}=0,
\qquad 1\leq i<j\leq k,
\]
to $\mathcal R_A$.
Denote the resulting set of relations by
$\mathcal R_M$.

For $1\leq j\leq k$, put
\[
c_j=[x,x_j]D\in M.
\]
We represent each $c_j$ by an $R$-module
coordinate vector with respect to the
generators $a_{ij}$.
If $q$ has finite order $o$, also put
\[
c_0=x^oD\in M.
\]

\begin{lemma}\label{lem:Cpresentation}
Fix $k\geq 1$.
Given a $k$-generated metabelian group $G$
defined by \eqref{eq:G} and an element
$x\in G$ represented by a word, the
$R$-module $C$ has the following presentation:
\[
C=
\begin{cases}
\left\langle
a_{ij}\ (i<j),\,e
\ \middle|\
\mathcal R_M,\;
(q_j-1)e=c_j\ (1\leq j\leq k)
\right\rangle_R,
& \text{if } |q|=\infty,\\[2mm]
\left\langle
a_{ij}\ (i<j),\,e
\ \middle|\
\mathcal R_M,\;
(q_j-1)e=c_j\ (1\leq j\leq k),\;
oe=c_0
\right\rangle_R,
& \text{if } |q|=o<\infty.
\end{cases}
\]
Here $e$ represents $xD$, and each $a_{ij}$
represents its image in $M$.

This presentation can be constructed in
deterministic polynomial time in $N+|x|$,
where $N=m+\sum_{i=1}^m |r_i|$,
with degree depending only on $k$.
\end{lemma}

\begin{proof}
Let $\widetilde C$ be the $R$-module
defined by the displayed presentation.
All its defining relations hold in $C$,
so the assignments
\[
a_{ij}\mapsto a_{ij}D,
\qquad
e\mapsto xD
\]
induce a surjective $R$-module homomorphism
\[
\psi:\widetilde C\to C.
\]
The natural map $M\to\widetilde C$ is
injective, since its composition with
$\psi$ is the inclusion $M\hookrightarrow C$.
Thus we identify $M$ with its image in
$\widetilde C$.

Let $\varepsilon:R\to\Z$ be the
augmentation map.
Since $\ker\varepsilon$ is generated
by $q_1-1,\ldots,q_k-1$, the defining
relations imply that
\[
re-\varepsilon(r)e\in M
\qquad\text{for every }r\in R.
\]
Consequently, every element of
$\widetilde C$ can be written as
\[
ze+m,\qquad z\in\Z,\quad m\in M.
\]

Suppose that $ze+m\in\ker\psi$.
Taking images in $Q=G/A$, we obtain
$q^z=1$.
If $q$ has infinite order, then $z=0$.
If $q$ has finite order $o$, then
$o\mid z$, and the relation $oe=c_0$
reduces $ze+m$ to an element of $M$.
In either case, the injectivity of
$\psi$ on $M$ implies that
$ze+m=0$.
Hence $\psi$ is an isomorphism.

We now establish the complexity bound.
By Proposition~\ref{prop:A-presentation},
an $R$-module presentation of $A$
can be constructed in polynomial time.
A presentation of $M=A/(q-1)A$ is
then obtained by adjoining the relations
\[
(q-1)a_{ij}=0,\qquad i<j.
\]
Their coefficients have polynomial-size
Laurent polynomial representatives.

Recall that
\[
Q\cong\Z^k/\operatorname{im}V,
\]
where $V$ is the matrix of exponent-sum
vectors of the defining relations.
Using Smith normal form, we can determine
whether $q$ has finite order and, if so,
compute its order $o$ in polynomial time.
The entries of $V$ have absolute value
at most $N$, and every nonzero minor
of $V$ has absolute value at most
$k!N^k$.

The order of every torsion element of $Q$
divides the largest nonzero invariant factor
of $V$. This invariant factor is bounded
by the greatest common divisor of the
nonzero maximal minors of $V$, and hence
by $k!N^k$.
Consequently,
\[
o\leq k!N^k.
\]
Thus, for fixed $k$, the order $o$
is polynomially bounded in $N$.

By Proposition~\ref{pr:collection},
the coordinate vectors of
\[
c_j=[x,x_j]D,\qquad 1\leq j\leq k,
\]
can be computed in polynomial time.

If $q$ has finite order $o$, then
$x^o\in A$. Since $o$ is polynomially
bounded, the word $x^o$ has polynomial
length in the combined input size.
Therefore Proposition~\ref{pr:collection}
also computes the coordinate vector
of $x^o$ in polynomial time.
Its image in $M$ gives the required
vector $c_0=x^oD$.

Finally, the presentation of $C$ is obtained
from that of $M$ by adding one generator
and at most $k+1$ relations.
The coefficients of these relations admit
Laurent polynomial representatives of
polynomial size.
Hence the presentation can be constructed
in deterministic polynomial time,
with the degree of the polynomial
depending only on $k$.
\end{proof}

\subsection{Polynomial-time reduction to a module orbit membership problem}

The next proposition is the main result of this section.

\begin{proposition}\label{prop:conj-to-orbit}
Fix $k\geq 1$.
There is a deterministic polynomial-time algorithm
that, given a $k$-generated metabelian group
$G$ defined by \eqref{eq:G} and two words
$x,y$ over its generators,
either determines that $x$ and $y$ are not conjugate,
or constructs
\begin{itemize}
\item
a finite $R$-module presentation of a module $C$
with $\binom{k}{2}+1$ generators;
\item
two elements $e,f\in C$, specified by $R$-module
coordinate vectors,
\end{itemize}
such that
\[
x\sim_G y
\quad\Longleftrightarrow\quad
f=e\cdot q
\quad\text{for some }q\in Q.
\]
\end{proposition}

\begin{proof}
Compute the images of $x$ and $y$ in $Q$.
If they are different, then $x$ and $y$ are not conjugate.
Otherwise, put
\[
q=\varphi(x)=\varphi(y),\qquad
D=(q-1)G',\qquad M=G'/D,
\]
and
\[
C=\langle M,xD\rangle\leq G/D.
\]

By Lemma~\ref{le:C-module}(b), $C$ is an
$R$-module, and by Proposition~\ref{prop:noskov},
conjugacy of $x$ and $y$ is equivalent to
the orbit condition for $xD$ and $yD$.
By Lemma~\ref{lem:Cpresentation}, a finite
$R$-module presentation of $C$ can be constructed
in polynomial time, with a distinguished generator
$e$ representing $xD$.
Since $x^{-1}y\in A$, Proposition~\ref{pr:collection}
computes an $R$-module coordinate vector $\alpha$
of $x^{-1}y$ in polynomial time.
Let $\bar\alpha$ be its image in $M$.
Then
\[
f=e+\bar\alpha
\]
represents $yD$ by Lemma~\ref{le:C-module}(c).

Thus
\[
x\sim_G y
\quad\Longleftrightarrow\quad
f=e\cdot s
\quad\text{for some }s\in Q.
\]
Since $k$ is fixed, the number of module
generators and the number of variables in
$R$ are bounded independently of the
defining relations of $G$.
All constructions take deterministic
polynomial time in the combined input size,
with degree depending only on $k$.
\end{proof}

\section{Reduction of orbit membership to relation-lattice membership}
\label{se:affine}

In this section, we reduce the module orbit problem
obtained in Section~\ref{sec:group} to subgroup
membership in the unit group of a commutative
affine $\Z$-algebra.
By Corollary~\ref{cor:affine-membership},
this problem can be solved in bounded-error
quantum polynomial time.

Recall that $Q=G/G'$ is generated by the images of
$x_1,\ldots,x_k$ and that $R=\Z Q$.
By \eqref{eq:R-presentation}, the ring $R$ admits
a presentation as a quotient of
\[
\Z[X_1,Y_1,\ldots,X_k,Y_k].
\]
Since $k$ is fixed, the number $2k$ of
polynomial variables is fixed, independently
of the defining relations of $G$.

We first reduce the orbit-membership problem
\[
f=e\cdot q
\qquad\text{for some }q\in Q,
\]
where $e,f$ belong to a finitely presented
$R$-module $C$, to subgroup membership in
the unit group of a suitable quotient ring.

Every element $q\in Q$ is a unit in $R=\Z Q$,
with inverse $q^{-1}$.
Its image in any quotient ring $S=R/I$
is therefore also a unit.
Consequently, the natural map $Q\to S$
takes values in $S^\times$ and defines
a group homomorphism
\[
Q\longrightarrow S^\times,
\qquad q\longmapsto q+I.
\]

\begin{proposition}\label{prop:orbit}
Let $C$ be a finitely presented $R=\Z Q$-module,
and let $e,f\in C$.
\begin{itemize}
\item[(a)]
If $e=0$, then
\[
f\in eQ \quad\Longleftrightarrow\quad f=0.
\]

\item[(b)]
Suppose that $e\ne0$. If $f\notin eR$, then
$f\notin eQ$. Otherwise, choose $h\in R$ such
that $f=eh$, and put
\[
I=\Ann_R(e),\qquad S=R/I.
\]
Then
\begin{equation}\label{eq:cyclic-orbit}
f\in eQ
\quad\Longleftrightarrow\quad
h+I\in\im\bigl(Q\longrightarrow S^\times\bigr).
\end{equation}
\end{itemize}
\end{proposition}

\begin{proof}
Part (a) is immediate, since the orbit of $0$
consists only of $0$.
For part (b), observe that $eQ\subseteq eR$.
Thus $f\notin eR$ implies $f\notin eQ$.
Suppose now that $f=eh$ for some $h\in R$.
For every $q\in Q$, we have
\[
f=eq
\quad\Longleftrightarrow\quad
e(h-q)=0
\quad\Longleftrightarrow\quad
h+I=q+I.
\]
This proves \eqref{eq:cyclic-orbit}.
\end{proof}

\begin{corollary}\label{co:orbit}
Fix $k\geq 1$.
There is a deterministic polynomial-time algorithm
that, given
\begin{itemize}
\item
a $k$-generated abelian group $Q$ defined by \eqref{eq:Q-presentation};
\item
a finite $R$-module presentation of a module $C$,
where $R=\Z Q$;
\item
two elements $e,f\in C$, specified by
$R$-module coordinate vectors,
\end{itemize}
performs the following computations:
\begin{itemize}
\item[(a)]
decides whether $f\in eR$ and, if so,
computes $h\in R$ such that $f=eh$;
\item[(b)]
computes generators of $I=\Ann_R(e)$;
\item[(c)]
constructs a quotient-ring presentation
of $S=R/I$;
\item[(d)]
if $h$ has been computed, decides whether
$h+I$ is a unit in $S$ and, if so,
computes its inverse.
\end{itemize}
\end{corollary}

\begin{proof}
The test $f\in eR$ is the membership problem for the
cyclic submodule $eR$ in $C$, and is therefore decided by
Corollary~\ref{cor:group-ring}(c).
If $f\in eR$, then Corollary~\ref{cor:group-ring}(c)
also computes $h\in R$ such that
\[
f=eh.
\]
Generators $g_1,\ldots,g_\ell$ of the ideal
$I=\Ann_R(e)$ are computed by
Corollary~\ref{cor:group-ring}(d).

Using \eqref{eq:R-presentation}, compute
polynomial representatives $\widetilde g_i$
of $g_i$ in $\Z[X_1,Y_1,\ldots,X_k,Y_k]$.
Then
\begin{equation}\label{eq:R/I-presentation}
S\cong
\frac{\Z[X_1,Y_1,\ldots,X_k,Y_k]}
{\left\langle
X_iY_i-1\ (1\le i\le k),\
\mono(r_j)-1\ (1\le j\le m),\
\widetilde g_1,\ldots,\widetilde g_\ell
\right\rangle}.
\end{equation}
This presentation of $S=R/I$ can be constructed in polynomial time.

Finally, Corollary~\ref{cor:inverse-mod-ideal} decides whether
$h+I$ is a unit in $S$ and, when it is, computes its inverse.
\end{proof}

\begin{proposition}[Reduction to affine membership]
\label{pr:orbit-to-affine}
Fix $k\geq 1$.
There is a deterministic polynomial-time algorithm
that, given
\begin{itemize}
\item
a $k$-generated metabelian group $G$
defined by \eqref{eq:G}, with
$Q=G/G'$ and $R=\Z Q$;
\item
a finite $R$-module presentation of $C$;
\item
two elements $e,f\in C$, specified by
$R$-module coordinate vectors,
\end{itemize}
either decides whether $f\in eQ$,
or constructs
\begin{itemize}
\item
a presentation \eqref{eq:R/I-presentation}
of a commutative affine $\Z$-algebra $S$;
\item
units $h,u_1,\ldots,u_k\in S^\times$,
specified by polynomial representatives
of themselves and their inverses,
\end{itemize}
such that
\[
f\in eQ
\quad\Longleftrightarrow\quad
h\in\langle u_1,\ldots,u_k\rangle.
\]
All constructed presentations and polynomial
representatives have polynomial size.
\end{proposition}

\begin{proof}
First decide whether $e=0$ by testing
membership in the zero submodule, using
Corollary~\ref{cor:group-ring}(c).
If $e=0$, then $eQ=\{0\}$, so membership
of $f$ in $eQ$ is decided by testing
whether $f=0$.

Suppose that $e\ne0$.
By Corollary~\ref{co:orbit}, we can decide
whether $f\in eR$.
If $f\notin eR$, then $f\notin eQ$.
Otherwise, compute $h_0\in R$ such that $f=eh_0$,
and compute generators of $I=\Ann_R(e)$.
Construct the quotient ring $S=R/I$ using \eqref{eq:R/I-presentation}.

By Corollary~\ref{cor:inverse-mod-ideal},
we can decide whether $h_0+I$ is a unit
in $S$ and compute its inverse when it is.
If $h_0+I$ is not a unit, then
Proposition~\ref{prop:orbit} implies
$f\notin eQ$.
Suppose therefore that $h=h_0+I$
is a unit in $S$.
Let $u_i=x_i+I\in S^\times$,
where $x_1,\ldots,x_k$ denote the
standard generators of $Q$.
In the polynomial presentation
\eqref{eq:R/I-presentation} of $S$,
the elements $u_i$ are represented
by $X_i$. Since $X_iY_i=1$ in $S$,
their inverses are represented by $Y_i$.

By Proposition~\ref{prop:orbit},
\[
f\in eQ
\quad\Longleftrightarrow\quad
h\in\im(Q\to S^\times).
\]
Since $Q$ is generated by
$x_1,\ldots,x_k$, its image in
$S^\times$ is precisely
$\langle u_1,\ldots,u_k\rangle$.
Therefore
\[
f=e\cdot q
\quad\text{for some }q\in Q
\quad\Longleftrightarrow\quad
h\in\langle u_1,\ldots,u_k\rangle.
\]

All preliminary tests and constructions
take deterministic polynomial time in the
combined input size, with degree depending
only on $k$, by Corollaries~\ref{cor:group-ring},
\ref{co:orbit}, and
\ref{cor:inverse-mod-ideal}.
The resulting presentations and polynomial
representatives have polynomial size.
\end{proof}

\section{Proof of Theorem~\ref{thm:group-main-uniform}}
\label{se:final-proof}

Fix $k\geq 1$. Consider a $k$-generated metabelian group
given by \eqref{eq:G} and
words $x,y$ over its generators.
Put
\[
A=G',\qquad Q=G/A,\qquad R=\Z Q,
\]
and
\[
N=m+\sum_{i=1}^m |r_i|,
\qquad L=N+|x|+|y|.
\]

By Proposition~\ref{prop:A-presentation},
an $R$-module presentation of $A$ can be
constructed in deterministic polynomial time
in $L$, with degree depending only on $k$.
By Proposition~\ref{prop:conj-to-orbit}
and Lemma~\ref{lem:Cpresentation}, we can
either determine in deterministic polynomial
time that $x$ and $y$ are not conjugate,
or construct a finitely presented
$R$-module $C$ and elements $e,f\in C$,
specified by $R$-module coordinate vectors,
such that
\[
x\sim_G y
\quad\Longleftrightarrow\quad
f=e\cdot q
\quad\text{for some }q\in Q.
\]
The presentation of $C$ and the coordinate
vectors of $e,f$ have polynomial size
and can be constructed in deterministic
polynomial time in $L$.

By Proposition~\ref{pr:orbit-to-affine}
and Corollary~\ref{co:orbit}, we can either
decide whether $f\in eQ$ in deterministic
polynomial time, or construct a commutative
affine $\Z$-algebra $S$ and units
$h,u_1,\ldots,u_k\in S^\times$ such that
\[
f\in eQ
\quad\Longleftrightarrow\quad
h\in\langle u_1,\ldots,u_k\rangle.
\]
The algebra $S$ is presented as a quotient
of a polynomial ring in $2k$ variables,
and the units and their inverses are
specified by polynomial representatives.
All these data have polynomial size
and can be constructed in deterministic
polynomial time in $L$.

Since $k$ is fixed, both the number $2k$
of polynomial variables and the number
$k$ of subgroup generators are fixed,
independently of the defining relations
of $G$.
Therefore, Corollary~\ref{cor:affine-membership}
decides whether
\[
h\in\langle u_1,\ldots,u_k\rangle
\]
in bounded-error quantum polynomial time
in $L$.

Combining the deterministic polynomial-time
reductions with this quantum membership
algorithm gives a bounded-error quantum
polynomial-time algorithm for deciding
whether $x$ and $y$ are conjugate in $G$,
with the degree of the polynomial depending
only on $k$.
\hfill$\square$

\bibliographystyle{plain}
\bibliography{bibliography}
\end{document}